\documentclass[12pt,oneside]{article}
\usepackage{amsmath,amssymb,amsfonts,amsthm}
\usepackage{color}

\newcommand{\T}{{\cal T}}

\newcommand{\Real}{\mathbb R}

\newcommand{\To}{\longrightarrow}

\newcommand{\prof}{\noindent \textit{\textbf{Proof.\:\:}}}
\newcommand{\tm}{\T M}
\newcommand{\p}{\pi^{-1}(TM)}
\newcommand {\cp}{\mathfrak{X}(\pi (M))}

\newcommand {\cpp}{\mathfrak{X}(\T M)}
\def\o#1{\overline{#1}}

\def\pa{\partial}
\def\paa{\dot{\partial}}

\def\Section#1{\vspace{30truept}\addtocounter{section}{1}\setcounter{thm}{0}\setcounter{equation}{0}
{\noindent\Large\bf\arabic{section}.~~#1}\par \vspace{12pt}}

\newtheorem{thm}{Theorem}[section]
\newtheorem{cor}[thm]{Corollary}
\newtheorem{lem}[thm]{Lemma}
\newtheorem{prop}[thm]{Proposition}
\newtheorem{defn}[thm]{Definition}

\newtheorem{rem}[thm]{Remark}

\numberwithin{equation}{section}

\begin{document}
\title{{\bf A GLOBAL APPROACH TO THE THEORY OF \\SPECIAL FINSLER MANIFOLDS} }
\author{{\bf Nabil L. Youssef$^{\dag}$, S. H. Abed$^{\dag}$ and A. Soleiman$^{\ddag}$}}
\date{}

\maketitle                     
\vspace{-1.15cm}
\begin{center}
{$^{\dag}$Department of Mathematics, Faculty of Science,\\ Cairo
University, Giza, Egypt.}
\end{center}
\vspace{-0.8cm}
\begin{center}
nyoussef@frcu.eun.eg,\ \ sabed@frcu.eun.eg
\end{center}
\vspace{-0.7cm}
\begin{center}
and
\end{center}
\vspace{-0.7cm}
\begin{center}
{$^{\ddag}$Department of Mathematics, Faculty of Science,\\ Benha
University, Benha, Egypt.}
\end{center}
\vspace{-0.8cm}
\begin{center}
soleiman@mailer.eun.eg
\end{center}

\smallskip
\begin{center}
{\textbf{Dedicated to the memory of Prof. Dr. A. TAMIM}}
\end{center}

\vspace{1cm} \maketitle
\smallskip

\noindent{\bf Abstract.}  The aim of the present paper is to provide
a \emph{global} presentation of the theory of special Finsler
manifolds. We introduce and investigate\emph{ globally} (or
intrinsically, free from local coordinates) many of the most
important and most commonly used special Finsler manifolds\,:
locally Minkowskian, Berwald, Landesberg, general Landesberg,
$P$-reducible, $C$-reducible, semi-$C$-reducible,
quasi-$C$-reducible, $P^{*}$-Finsler, $C^{h}$-recurrent,
$C^{v}$-recurrent, $C^{0}$-recurrent, $S^{v}$-recurrent,
$S^{v}$-recurrent of the second order, $C_{2}$-like, $S_{3}$-like,
$S_{4}$-like, $P_{2}$-like, $R_{3}$-like, $P$-symmetric,
$h$-isotropic, of scalar curvature, of constant curvature, of
$p$-scalar curvature, of $s$-$ps$-curvature.
\par The global definitions of these special Finsler manifolds are
introduced. Various relationships between the different types of the
considered special Finsler manifolds are found. Many local results,
known in the literature, are proved globally and several new results
are obtained. As a by-product, interesting identities and properties
concerning the torsion tensor fields and the curvature tensor fields
are deduced.
\par Although our investigation is entirely global, we
provide; for comparison reasons, an appendix presenting a local
counterpart of our global approach and the {\it{local}} definitions
of the special Finsler spaces considered. \footnote{ArXiv Number:
0704.0053}

\bigskip
\medskip\noindent{\bf Keywords and phrases.\/}\,  Berwald, Landesberg,
$P$-reducible, $C$-reducible, Semi-$C$-reducible,
Quasi-$C$-reducible, $P^{*}$-Finsler, $C^{h}$-recurrent,
$C^{v}$-recurrent, $S^{v}$-recurrent, $C_{2}$-like, $S_{3}$-like,
$S_{4}$-like, $P_{2}$-like, $R_{3}$-like, $P$-symmetric,
$h$-isotropic, Of scalar curvature, Of constant curvature, Of
$p$-scalar curvature, Of $s$-$ps$-curvature.

\bigskip
\medskip\noindent{\bf  2000 AMS Subject Classification.\/} 53C60,
53B40.
\newpage
\vspace{30truept}\centerline{\Large\bf{Introduction}}\vspace{12pt}
\par \par In Finsler geometry all
geometric objects  depend not only on positional coordinates, as in
Riemannian geometry, but also on directional arguments. In
Riemannian geometry there is a canonical linear connection on the
manifold $M$, while in Finsler geometry there is a corresponding
canonical linear connection, due to E. Cartan, which is not a
connection on $M$ but is a connection on $\,\pi^{-1}(TM) $,  the
pullback of the tangent bundle $TM$ by $\,\pi: \T M\longrightarrow
M$ (\emph{the pullback approach}). Moreover, in Riemannian geometry
 there is one curvature tensor and one torsion tensor  associated with a given
 linear connection on the manifold $M$, whereas  in Finsler geometry
 there are three curvature tensors  and five torsion tensors
 associated with a given linear connection on $\,\pi^{-1}(TM) $.
\par
Most of the special spaces in  Finsler geometry are derived from the
fact that the $\pi$-tensor fields (torsions and curvatures)
associated with the Cartan connection satisfy special forms.
Consequently, special spaces of Finsler geometry are more numerous
than those of Riemannian geometry. Special Finsler spaces are
investigated locally (using local coordinates) by many authors: M.
Matsumoto \cite{r34}, \cite{r32}, \cite{r2}, \cite{r29} and others
\cite{r10}, \cite{r31}, \cite{r75}, \cite{r76}. On the other hand,
the global (or intrinsic, free from local coordinates) investigation
of such spaces is very rare in the literature. Some considerable
contributions in this direction are due to A. Tamim \cite{r44},
\cite{r48}.
\par
In the present paper, we provide a \emph{global} presentation of the
theory of special Finsler manifolds. We introduce and
investigate\emph{ globally} many of the most important and most
commonly used special Finsler manifolds\,: locally Minkowskian,
Berwald, Landesberg, general Landesberg, $P$-reducible,
$C$-reducible, semi-$C$-reducible, quasi-$C$-reducible,
$P^{*}$-Finsler, $C^{h}$-recurrent, $C^{v}$-recurrent,
$C^{0}$-recurrent, $S^{v}$-recurrent, $S^{v}$-recurrent of the
second order, $C_{2}$-like, $S_{3}$-like, $S_{4}$-like,
$P_{2}$-like, $R_{3}$-like, $P$-symmetric, $h$-isotropic, of scalar
curvature, of constant curvature, of $p$-scalar curvature, of
$s$-$ps$-curvature.
\par
The paper consists of two parts, preceded by a preliminary section
$(\S 1)$, which provides a brief account of the basic concepts of
the pullback approach to Finsler geometry necessary to this work.
For more detail, the reader is referred to \cite{r57}, \cite{r58},
\cite{r61} and \cite{r44}.
 \par
In the first part $(\S 2)$, we introduce the global definitions of
the aforementioned special Finsler manifolds in such a way that,
when localized, they yield the usual local definitions current in
the literature (see the appendix). The definitions are arranged
according to the type of the defining property of the special
Finsler manifold concerned.
\par
In the second part $(\S 3)$, various relationships between the
different types of the considered special Finsler manifolds are
found. Many local results, known in the literature, are proved
globally and several new results are obtained. As a by-product of
some of the obtained results, interesting identities and properties
concerning the torsion tensor fields and the curvature tensor fields
are deduced, which in turn play a key role in obtaining other
results.
\par
Among the obtained results are: a characterization of Riemannian
manifolds, a characterization of $S^{v}$-recurrent manifolds, a
characterization of $P$-symmetric \linebreak manifolds, a
characterization of Berwald manifolds (in certain cases), the
equivalence of Landsberg and general Landsberg manifolds under
certain conditions, a classification of $h$-isotropic
$C^{h}$-recurrent manifolds and a presentation of different
conditions under which an $R_{3}$-like Finsler manifold becomes a
Finsler manifold of $s$-$ps$ curvature. The above results are just a
non-exhaustive sample of the global results obtained in this paper.
\par
It should finally be noted that some important results of
\cite{r75}, \cite{r84}, \cite{r69}, \cite{r65}, \cite{r31},
\cite{r68},...,etc. (obtained in local coordinates) are immediately
derived from the obtained global results (when localized).
\par
Although our investigation is entirely global, we conclude the paper
with an appendix presenting a local counterpart of our global
approach and the {\it{local}} definitions of the special Finsler
spaces considered. This is done to facilitate comparison and to make
the paper more self-contained.

\Section{Notation and Preliminaries}

In this section, we give a brief account of the basic concepts
 of the pullback\linebreak formalism of Finsler geometry necessary for this work. For more details
  refer to~\cite{r57},\,\cite{r58},\linebreak \,\cite{r61} and~\,\cite{r44}.
 We make the general
assumption that all geometric objects we consider are of class
$C^{\infty}$.
The following notations will be used throughout this paper:\\
 $M$: a real differentiable manifold of finite dimension $n$ and of
class $C^{\infty}$,\\
 $\mathfrak{F}(M)$: the $\Real$-algebra of differentiable functions
on $M$,\\
 $\mathfrak{X}(M)$: the $\mathfrak{F}(M)$-module of vector fields
on $M$,\\
$\pi_{M}:TM\longrightarrow M$: the tangent bundle of $M$,\\
$\pi: \T M\longrightarrow M$: the subbundle of nonzero vectors
tangent to $M$,\\
$V(TM)$: the vertical subbundle of the bundle $TTM$,\\
 $P:\pi^{-1}(TM)\longrightarrow \T M$ : the pullback of the
tangent bundle $TM$ by $\pi$,\\
$P^*:\pi^{-1}(T^{*}M)\longrightarrow \T M$ : the pullback of the
cotangent bundle $T^{*}M$ by $\pi$,\\
 $\mathfrak{X}(\pi (M))$: the $\mathfrak{F}(TM)$-module of
differentiable sections of  $\pi^{-1}(T M)$.\\
\par Elements  of  $\mathfrak{X}(\pi (M))$ will be called
$\pi$-vector fields and will be denoted by barred letters
$\overline{X} $. Tensor fields on $\pi^{-1}(TM)$ will be called
$\pi$-tensor fields. The fundamental $\pi$-vector field is the
$\pi$-vector field $\overline{\eta}$ defined by
$\overline{\eta}(u)=(u,u)$ for all $u\in \T M$. The lift to
$\pi^{-1}(TM)$ of a vector field $X$ on $M$ is the $\pi$-vector
field $\overline{X}$ defined by $\overline{X}(u)=(u,X(\pi (u)))$.
The lift to $\pi^{-1}(TM)$ of a $1$-form $\omega$ on $M$ is the
$\pi$-form $\overline{\omega}$ defined by
$\overline{\omega}(u)=(u,\omega(\pi (u)))$.

 The tangent bundle $T(\T M)$
is related to the pullback bundle $\pi^{-1}(TM)$ by the short exact
sequence \vspace{-0.4cm}
$$0\longrightarrow
 \pi^{-1}(TM)\stackrel{\gamma}\longrightarrow T(\T M)\stackrel{\rho}\longrightarrow
\pi^{-1}(TM)\longrightarrow 0 ,\vspace{-0.2cm}$$
 where the bundle morphisms $\rho$ and $\gamma$ are defined respectively by
$\rho = (\pi_{\T M},d\pi)$ and $\gamma (u,v)=j_{u}(v)$, where
$j_{u}$  is the natural isomorphism $j_{u}:T_{\pi_{M}(v)}M
\longrightarrow T_{u}(T_{\pi_{M}(v)}M)$.

 Let $\nabla$ be  a linear
connection (or simply a connection) in the pullback bundle
$\pi^{-1}(TM)$.
 We associate to
$\nabla$ the map\vspace{-0.2cm}
$$K:T \T M\longrightarrow \pi^{-1}(TM):X\longmapsto \nabla_X \overline{\eta}
,\vspace{-0.2cm}$$ called the connection (or the deflection) map of
$\nabla$. A tangent vector $X\in T_u (\T M)$ is said to be
horizontal if $K(X)=0$ . The vector space $H_u (\T M)= \{ X \in T_u
(\T M) : K(X)=0 \}$ of the horizontal vectors
 at $u \in  \T M$ is called the horizontal space to $M$ at $u$  .
   The connection $\nabla$ is said to be regular if
 $$T_u (\T M)=V_u (\T M)\oplus H_u (\T M) \qquad \forall u\in \T M .$$
  \par If $M$ is endowed with a regular connection, then the vector bundle
   maps
\begin{eqnarray*}
 \gamma &:& \pi^{-1}(T M)  \To V(\T M), \\
   \rho |_{H(\T M)}&:&H(\T M) \To \pi^{-1}(TM), \\
   K |_{V(\T M)}&:&V(\T M) \To \pi^{-1}(T M)
\end{eqnarray*}
 are vector bundle isomorphisms.
   Let us denote
 $\beta=(\rho |_{H(\T M)})^{-1}$,
then \vspace{-0.2cm}
   \begin{align}\label{fh1}
    \rho  o  \beta = id_{\pi^{-1} (\T M)}, \quad  \quad
       \beta o \rho =\left\{
                                \begin{array}{ll}
                                          id_{H(\T M)} & \hbox{on   H(\T M)} \\
                                         0 & \hbox{on    V(\T M)}
                                       \end{array}
                                     \right.\vspace{-0.2cm}
\end{align}

For a regular connection $\nabla$ we define two covariant
derivatives $\stackrel{1}\nabla$ and $\stackrel{2}\nabla$ as
follows: For every
 vector (1)$\pi$-form $A$, we have\vspace{-0.2cm}
 $$ (\stackrel{1}\nabla A)(\o X, \o Y):= (\nabla_{\beta \o X} A)( \o Y)\,, \qquad\qquad
  (\stackrel{2}\nabla A)(\o X, \o Y):= (\nabla_{\gamma \o X} A)( \o Y).$$

 The classical  torsion tensor $\textbf{T}$  of the connection
$\nabla$ is defined by
$$\textbf{T}(X,Y)=\nabla_X \rho Y-\nabla_Y\rho X -\rho [X,Y] \quad
\forall\,X,Y\in \mathfrak{X} (\T M).$$ The horizontal ((h)h-) and
mixed ((h)hv-) torsion tensors, denoted respectively by $Q $ and $ T
$, are defined by \vspace{-0.2cm}
$$Q (\overline{X},\overline{Y})=\textbf{T}(\beta \overline{X}\beta \overline{Y}),
\, \,\, T(\overline{X},\overline{Y})=\textbf{T}(\gamma
\overline{X},\beta \overline{Y}) \quad \forall \,
\overline{X},\overline{Y}\in\mathfrak{X} (\pi (M)).\vspace{-0.2cm}$$
\par The classical curvature tensor  $\textbf{K}$ of the connection
$\nabla$ is defined by
 $$ \textbf{K}(X,Y)\rho Z=-\nabla_X \nabla_Y \rho Z+\nabla_Y \nabla_X \rho Z+\nabla_{[X,Y]}\rho Z
  \quad \forall\, X,Y, Z \in \mathfrak{X} (\T M).$$
The horizontal (h-), mixed (hv-) and vertical (v-) curvature
tensors, denoted respectively by $R$, $P$ and $S$, are defined by
$$R(\overline{X},\overline{Y})\o Z=\textbf{K}(\beta
\overline{X}\beta \overline{Y})\o Z,\quad
P(\overline{X},\overline{Y})\o Z=\textbf{K}(\beta
\overline{X},\gamma \overline{Y})\o Z,\quad
S(\overline{X},\overline{Y})\o Z=\textbf{K}(\gamma
\overline{X},\gamma \overline{Y})\o Z.$$ We also have the
 (v)h-, (v)hv- and (v)v-torsion tensors,
denoted respectively by $\widehat{R}$, $\widehat{P}$ and
$\widehat{S}$, defined by
$$\widehat{R}(\overline{X},\overline{Y})={R}(\overline{X},\overline{Y})\o \eta,\quad
\widehat{P}(\overline{X},\overline{Y})={P}(\overline{X},\overline{Y})\o
\eta,\quad
\widehat{S}(\overline{X},\overline{Y})={S}(\overline{X},\overline{Y})\o
\eta.$$
\begin{thm} \label{th.1}{\em{\cite{r48}}} Let $(M,L)$ be a Finsler manifold. There exists a
unique regular connection $\nabla$ in $\pi^{-1}(TM)$ such that
\begin{description}
  \item[(a)]  $\nabla$ is  metric\,{\em:} $\nabla g=0$,

  \item[(b)]   The horizontal torsion of $\nabla$ vanishes\,{\em:} $Q=0
  $,
  \item[(c)]  The mixed torsion $T$ of $\nabla$ satisfies \,
  $g(T(\overline{X},\overline{Y}), \overline{Z})=g(T(\overline{X},\overline{Z}),\overline{Y})$.
\end{description}
\end{thm}
Such a connection is called the Cartan connection associated to the
Finsler manifold $(M,L)$.

\vspace{6pt} One can show that the torsion $T$ of the Cartan
connection has the property that $T(\overline{X},\overline{\eta})=0$
for all $\overline{X} \in \mathfrak{X} (\pi (M))$ and associated to
$T$ we have:

\begin{defn}{\em{\cite{r48}}} Let  $\nabla$ be the Cartan  connection associated to $(M,L)$.
The torsion tensor field  $T$ of the connection $\nabla$ induces a
$\pi$-tensor field of type $(0,3)$, called the Cartan tensor and
denoted again $T$, defined by{\,\em{:}}
$$T(\overline{X},\overline{Y},\overline{Z})=g(T(\overline{X},\overline{Y}),\overline{Z}),
\quad \text {for all} \quad
\overline{X},\overline{Y},\overline{Z}\in \mathfrak{X} (\T M).$$ It
also induces a $\pi$-form $C$, called the contracted torsion,
defined by{\,\em{:}}
$$C(\overline{X}):= Tr\{\overline{Y} \longmapsto
T(\overline{X},\overline{Y})\},
 \quad \text {for all} \quad \overline{X}\in \mathfrak{X} (\T M).$$
\end{defn}

\begin{defn} {\em{\cite{r48}}} With respect to the Cartan  connection $\nabla$ associated to
$(M,L)$, we have
 \begin{description}
    \item[--]The horizontal and vertical Ricci tensors $Ric^h$ and
    $Ric^v$  are defined respectively by:
 $$Ric^h(\overline{X},\overline{Y}):= Tr\{ \overline{Z} \longmapsto
R(\overline{X},\overline{Z})\overline{Y}\}, \quad \text {for all}
\quad \overline{X},\overline{Y}\in \mathfrak{X} (\T M),$$
 $$Ric^v(\overline{X},\overline{Y}):= Tr\{ \overline{Z} \longmapsto
S(\overline{X},\overline{Z})\overline{Y}\}, \quad \text {for all}
\quad \overline{X},\overline{Y}\in \mathfrak{X} (\T M).$$

   \item[--]The horizontal and vertical Ricci  maps $Ric_0^h$ and $Ric_0^v$ are
defined respectively by:
$$g(Ric_0^h(\overline{X}),\overline{Y}):=Ric^h(\overline{X},\overline{Y}),
\quad \text {for all} \quad \overline{X},\overline{Y}\in
\mathfrak{X} (\T M),$$
$$g(Ric_0^v(\overline{X}),\overline{Y}):=Ric^v(\overline{X},\overline{Y}),
\quad \text {for all} \quad \overline{X},\overline{Y}\in
\mathfrak{X} (\T M).$$

    \item[--] The horizontal and vertical scalar
curvatures $Sc^h$ , $Sc^v$ are defined respectively
by:\vspace{-0.2cm}
\begin{equation*}
     Sc^h:=
\text {Tr}( Ric_0^h), \qquad \qquad  Sc^v:= \text {Tr}  (Ric_0^v
),\vspace{-0.3cm}
\end{equation*}
 \end{description}
where $R$ and $S$ are respectively the horizontal and vertical
curvature tensors of $\nabla$.
\end{defn}


\begin{prop}{\em{\cite{r27}}} Let $(M,L)$ be a Finsler manifold. The vector field $G$
 determined by $i_{G}\Omega =-dE$ is a spray,
 called the canonical spray associated to the energy $E$, where
 $E:=\frac{1}{2}L^{2}$ and $\Omega:=dd_{J}E$.
 \end{prop}
 One can show, in this case, that $G=\beta o \overline{\eta}$, and
 $G$ is thus horizontal with respect to the Cartan connection $\nabla$.

\begin{thm} {\em{\cite{r49}}} \label{th.1a} Let $(M,L)$ be a Finsler manifold. There exists  a
unique regular connection $D$ in $\pi^{-1}(TM)$ such that
\begin{description}
  \item[(a)]  $D$ is torsion free,
  \item[(b)]  The canonical spray $G= \beta o \overline{\eta}$ is
  horizontal with respect to $D$,
  \item[(c)] The (v)hv-torsion tensor $\widehat{P}$ of $D$ vanishes.
\end{description}
\end{thm}
Such a connection is called the Berwald connection associated to the
Finsler manifold $(M,L)$.



\newpage
\Section{Special Finsler spaces}

In this section, we introduce the global definitions of the most
important and commonly used special Finsler spaces in such a way
that, when localized, they yield the usual local definitions
existing in the literature (see the Appendix). Here we simply set
the definitions, postponing investigation of the mutual
relationships between these special Finsler spaces to the next
section. The definitions are arranged according to the type of
defining property of the special Finsler space concerned.
\par
Throughout the paper, $g$, $\widehat{g}$, $\nabla$ and $D$ denote
respectively the Finsler metric in $\pi^{-1}(TM)$, the induced
metric in $\pi^{-1}(T^{*}M)$, the Cartan connection and the Berwald
connection associated to a given Finsler manifold $(M,L)$. Also, $T$
denotes the torsion tensor of the Cartan connection (or the Cartan
tensor) and $R$, $P$ and $S$ denote respectively the horizontal
curvature, the mixed curvature and the vertical curvature of the
Cartan connection.

\begin{defn}\label{def.1a} A Finsler manifold $(M,L)$ is\,{\em:}
\begin{description}
  \item[(a)]  Riemannian  if the metric tensor $g(x,y)$ is independent
  of $y$ or, equivalently, if
  $$T(\overline{X},\overline{Y})=0,\,\,\,
    \text{for all}\,\,\,\overline{X}, \overline{Y}\in \mathfrak{X}(\pi(M)).$$

  \item[(b)] locally Minkowskian  if the metric tensor $g(x,y)$ is independent
  of $x$ or,  equivalently, if $$\nabla_{\beta \overline{X}}\,T
  =0\,\,\, \text{and}\,\,\,  R=0.$$
      \end{description}
\end{defn}

\begin{defn}\label{def.1aa} A Finsler manifold $(M,L)$ is said to be\,{\em:}
\begin{description}

 \item[(a)]  Berwald {\em{\cite{r44}}} if the torsion tensor $T$ is horizontally
  parallel. That is,
  $$\nabla_{\beta \overline{X}}\,T =0.$$

\item[(b)]  $C^h$-recurrent  if the torsion tensor
 $T$  satisfies the condition
  $$\nabla_{\beta \overline{X}}\,T
 =\lambda_{o} (\overline{X})\,T,$$
 \noindent where $\lambda_{o}$ is a $\pi$-form  of order one.
\item[(c)]
  $P^{*}$-Finsler manifold
  if the
$\pi$-tensor field  $\nabla_{\beta \overline{\eta}}T$
  is expressed in the form
  $$\nabla_{\beta \overline{\eta}}\,T= \lambda (x,y)\,T,$$
  where $\lambda (x,y)=
  \frac{\widehat{g}(\nabla_{\beta\overline{\eta}}\,C,C)}{C^2}=
  \frac{{g}(\nabla_{\beta\overline{\eta}} \o C, \o C)}{C^2}$  and
  $C^2:=\widehat{g}(C,C)=C(\overline{C})\neq 0$; $\overline{C}$ being the
  $\pi$-vector field defined by
  $g(\overline{C},\overline{X})=C(\overline{X})$.

\end{description}
\end{defn}

\begin{defn}A Finsler manifold $(M,L)$ is said to be{\em{\,\!:}}
\begin{description}

    \item[(a)]  $C^v$-recurrent  if the torsion tensor
 $T$  satisfies the condition\\
  $(\nabla_{\gamma \overline{X}}T)(\overline{Y},\overline{Z})
 =\lambda_{o} (\overline{X}) T(\overline{Y},\overline{Z}).$

    \item[(b)]  $C^0$-recurrent  if the torsion tensor
 $T$  satisfies the condition\\
  $(D_{\gamma \overline{X}}T)(\overline{Y},\overline{Z})
 =\lambda_{o} (\overline{X}) T(\overline{Y},\overline{Z}).$

\end{description}
\end{defn}

\begin{defn}{\em{\cite{r48}}}\label{3.def.1} A Finsler manifold $(M,L)$ is said to
be\,{\em:}
\begin{description}

 \item[(a)]  semi-$C$-reducible if $dim M \geq 3$ and the
  Cartan tensor $T$ has the form
  \begin{equation*}
  \begin{split}
    T(\overline{X},\overline{Y},\overline{Z})= & \frac{\mu}{n+1} \{\hbar(\overline{X}
  ,\overline{Y})C(\overline{Z})+\hbar(\overline{Y}
  ,\overline{Z})C(\overline{X})+\hbar(\overline{Z},\overline{X}) C(\overline{Y})\}+ \\
      &+ \frac{\tau}{C^{2}}C(\overline{X}) C(\overline{Y})
  C(\overline{Z}),
  \end{split}
  \end{equation*}
  \end{description}
where $\mu$ and $\tau$ are scalar functions satisfying $\mu
+\tau=1$,
   $\hbar= g-\ell \otimes \ell$ and
   $\ell(\overline{X}):=L^{-1}g(\overline{X},\overline{\eta})$.
\begin{description}
\item[(b)] $C$-reducible   if $dim M \geq 3$ and the
  Cartan tensor $T$ has the form
  \begin{equation*}
      T(\overline{X},\overline{Y},\overline{Z})=  \frac{1}{n+1} \{\hbar(\overline{X}
  ,\overline{Y})C(\overline{Z})+\hbar(\overline{Y}
  ,\overline{Z})C(\overline{X})+\hbar(\overline{Z},\overline{X})
  C(\overline{Y})\}.
   \end{equation*}

\item[(c)] $C_{2}$-like   if $dim M \geq 2$ and  the
  Cartan tensor $T$ has the form
  \begin{equation*}
    T(\overline{X},\overline{Y},\overline{Z})= \frac{1}{C^{2}}C(\overline{X}) C(\overline{Y})
  C(\overline{Z}).
  \end{equation*}

  \end{description}
\end{defn}

\begin{defn} \label{3.def.2} A Finsler manifold $(M,L)$, where  $dim M \geq 3$, is said to be
  quasi-$C$-reducible
if  the Cartan tensor $T$ is written as\,{\em:}
  $$T(\overline{X},\overline{Y},\overline{Z})= A(\overline{X}
  ,\overline{Y})C(\overline{Z})+A(\overline{Y}
  ,\overline{Z})C(\overline{X})+A(\overline{Z}
  ,\overline{X})C(\overline{Y}),$$
where $A$ is a symmetric indicatory  $(2)\,\pi$-form
{\em$(A(\overline{X},\overline{\eta})=0$ for all $\overline{X}$)}.
\end{defn}

\begin{defn} {\em{\cite{r48}}}\label{3.def.3} A Finsler manifold $(M,L)$ is said to be\,{\em:}
\begin{description}
\item[(a)]  $S_{3}$-like if $dim(M)\geq 4$ and the vertical curvature tensor
$S(\overline{X},\overline{Y},\overline{Z},\overline{W})$\\
$:=g(S(\overline{X},\overline{Y})\overline{Z},\overline{W})$ has the
form\,{\em:}
$$S(\overline{X},\overline{Y},\overline{Z},\overline{W})=
\frac{Sc^{v}}{(n-1)(n-2)} \{
\hbar(\overline{X},\overline{Z})\hbar(\overline{Y},\overline{W})-\hbar(\overline{X},\overline{W})
\hbar(\overline{Y},\overline{Z}) \}.$$

\item[(b)]  $S_{4}$-like  if $dim(M)\geq 5$ and the vertical curvature tensor
$S(\overline{X},\overline{Y},\overline{Z},\overline{W})$ has the
form\,{\em:}
\begin{equation}\label{h}
   \begin{split}
    S(\overline{X},\overline{Y},\overline{Z},\overline{W})=
    &
   \hbar(\overline{X},\overline{Z})\textbf{F}(\overline{Y},\overline{W})
   -\hbar(\overline{Y},\overline{Z})\textbf{F}(\overline{X},\overline{W})+ \\
       &+\hbar(\overline{Y},\overline{W})\textbf{F}(\overline{X},\overline{Z})
        - \hbar(\overline{X},\overline{W})\textbf{F}(\overline{Y},\overline{Z}),
   \end{split}
\end{equation}
\end{description}
\noindent where  $\textbf{F}$ is the $(2)\pi$-form defined by
$\textbf{F}= {\displaystyle\frac{1}{n-3}\{Ric^v- \frac{Sc^v\,
\hbar}{2(n-2)}\}}$.
\end{defn}

\begin{defn}\label{def.a1} A Finsler manifold $(M,L)$ is said to be\,{\em:}
\begin{description}

 \item[(a)]  $S^v$-recurrent  if the $v$-curvature  tensor
 $S$  satisfies the condition\\
  $$(\nabla_{\gamma \overline{X}}S)(\overline{Y},\overline{Z},\overline{W})
  =\lambda(\overline{X})S(\overline{Y},\overline{Z})\overline{W},$$
 where $\lambda$ is a $\pi$-form of order one.

\item[(b)]  $S^v$-recurrent  of the second order  if the $v$-curvature  tensor
 $S$  satisfies the condition
  $$(\stackrel{2}\nabla \stackrel{2} \nabla S)(\o Y, \o X, \overline{Z},\overline{W},\overline{U})
 =\Theta (\overline{X},\overline{Y})S(\overline{Z},\overline{W})\overline{U},$$
 where $\Theta$ is a $\pi$-form of order two.

\end{description}
\end{defn}

\newpage
\begin{defn} {\em{\cite{r44}}} \label{def.2a} A Finsler manifold $(M,L)$ is said to
be\,{\em:}
\begin{description}
\item[(a)] a Landsberg manifold  if
$$\widehat{P}(\overline{X},\overline{Y})=P(\overline{X},\overline{Y})\overline{\eta}=0\,\,\,
 \forall\,  \overline{X}, \overline{Y}\in\cp,
 \text {\, or equivalently \, }
 \nabla_{\beta \overline{\eta}}\,T =0.$$

 \item[(b)] a general Landsberg manifold if $$Tr\{
  \overline{Y} \To \widehat{P}(\overline{X},\overline{Y})\}=0\,\,\,
 \forall\,  \overline{X}, \in\cp,
  \text {\, or equivalently \, } \nabla _{\beta \overline{\eta}}\ C =0.$$
  \end{description}
\end{defn}


\begin{defn} \label{p-sy} A Finsler manifold $(M,L)$ is said to be
$P$-symmetric if the mixed curvature  tensor $P$
satisfies\vspace{-0.2cm}
$$P(\overline{X},\overline{Y})\overline{Z}
 =P(\overline{Y},\overline{X})\overline{Z}, \ \ \forall\  \o X, \o Y, \o Z\in\cp  .$$
\end{defn}

\begin{defn}\label{def.p2like}  A Finsler manifold $(M,L)$, where  $dim M\geq 3$,
is said to be $P_{2}$-like if the mixed curvature tensor $P$ has the
form\,{\em:}\vspace{-0.2cm}
$$P(\overline{X},\overline{Y},\overline{Z}, \o W)=
 \alpha(\overline{Z})T (\overline{X},\overline{Y},\o W)
 -\alpha(\overline{W})\, T(\overline{X},\o Y,\overline{Z}),\vspace{-0.2cm}$$ where $\alpha$ is
  a $(1)\,\pi$-form {\em{(}positively homogeneous of degree $0$)}.
  \end{defn}

\begin{defn} {\em{\cite{r48}}}\label{def.p1} A Finsler manifold $(M,L)$, where  $dim M\geq 3$,
is said to be $P$-reducible
  if  the $\pi$-tensor field
 $P(\overline{X},\overline{Y},\overline{Z}):
  =g(P(\overline{X},\overline{Y})\overline{\eta},\overline{Z})$ can be expressed
  in the form\,{\em:}
 $$P(\overline{X},\overline{Y},\overline{Z})=\delta(\overline{X})\hbar (\overline{Y},\overline{Z})
  +\delta(\overline{Y})\hbar (\overline{Z},\overline{X})+ \delta(\overline{Z})\hbar
  (\overline{X},\overline{Y}),$$ where $\delta$ is
  a $(1)\,\pi$-form satisfying   $\delta (\o \eta)=0.$
\end{defn}


\begin{defn}\label{def.2}{\em{\cite{r14}}} A Finsler manifold $(M,L)$, where  $dim M\geq 3$,
is said to be\linebreak $h$-isotropic if there exists a scalar
$k_{o}$ such that the horizontal curvature tensor $R$ has the form
$$R(\overline{X},\overline{Y})\overline{Z}=k_{o} \{g(\overline{Y},\overline{Z})
\overline{X}-g(\overline{X},\overline{Z})\overline{Y} \}.$$
\end{defn}

\begin{defn}\label{def.2}{\em{\cite{r14}}} A Finsler manifold $(M,L)$, where  $dim M\geq 3$, is said to be\,{\em:}
\begin{description}

 \item[(a)]  of scalar curvature  if there exists a  scalar function $k: \T M \To
 \Real$ such that the horizontal curvature tensor $R(\overline{X},\overline{Y},\overline{Z},\overline{W}):=
g(R(\overline{X},\overline{Y})\overline{Z},\overline{W})$ satisfies
the relation
 $$R(\overline{\eta},\overline{X},\overline{\eta},\overline{Y})=
  k L^{2} \hbar(\overline{X},\overline{Y}). $$

\item[(b)]  of constant curvature  if the function $k$
in {\em(a)} is constant.
\end{description}
\end{defn}


\begin{defn}\label{def.r3} A Finsler manifold $(M,L)$ is said to be $R_{3}$-like if
 $dim M\geq 4$ and the horizontal
curvature tensor
$R(\overline{X},\overline{Y},\overline{Z},\overline{W})$ is
expressed in the form\vspace{-0.2cm}
\begin{equation}\label{h}
   \begin{split}
    R(\overline{X},\overline{Y},\overline{Z},\overline{W})= &
       g(\overline{X},\overline{Z})F(\overline{Y},\overline{W})
   -g(\overline{Y},\overline{Z})F(\overline{X},\overline{W})+  \\
       &+ g(\overline{Y},\overline{W})F (\overline{X},\overline{Z})
       - g(\overline{X},\overline{W})F(\overline{Y},\overline{Z}) ,
   \end{split}
\end{equation}
\end{defn}\vspace{-0.3cm}
\noindent where $F$ is the $(2)\pi$-form defined by $ F={
\frac{1}{n-2}\{ Ric^h- \frac{Sc^h\, g}{2(n-1)}\}}$.
\newpage


\Section{Relationships between different types of special\\
Finsler spaces}

This section is devoted to global investigation of some mutual
relationships\linebreak  between the special Finsler spaces
introduced in the preceding section. Some consequences are also
drawn from these relationships. \vspace{3pt}
\par We start with some immediate consequences from the
definitions:\\
(a) A Locally  Minkowskian manifold is a Berwald  manifold.\\
(b) A Berwald manifold is a Landsberg  manifold.\\
(c) A  Landsberg manifold is a general Landsberg  manifold.\\
(d) A Berwald manifold is $C^{h}$-recurrent (resp. $P^{*}$-Finsler).\\
(e) A $P^{*}$-manifold is a Landsberg manifold.\\
(f) A $C$-reducible (resp. $C_{2}$-like) manifold  is semi-$C$-reducible.\\
(g) A semi-$C$-reducible manifold is quasi-$C$-reducible.\\
(h) A Finsler manifold of constant curvature is of scalar curvature.

\par
\vspace{7pt}
 The following two lemmas are useful for subsequent use.\vspace{-0.2cm}

\begin{lem}\label{le.p}{\em\cite{r48}} For every $\o X, \o Y\in \cp $, we
have{\em\,\!{:}}
$${ \textbf{\em(a)}}\ P(\o \eta, \o X)\o Y=0,\qquad\quad\
{\textbf{\em(b)}}\ P(\o X, \o \eta)\o Y=0,\qquad\quad\
{\textbf{\em(c)}}\ P(\o X, \o Y)\o \eta=(\nabla_{\beta \o \eta}T)(\o
X,\o Y).
$$
\end{lem}
\vspace{-0.3cm}
\begin{lem}\label{cor.i1} If $\phi$ is the  vector $\pi$-form  defined  by
\vspace{-0.2cm}
\begin{equation}\label{eq.i1}
   \phi(\o X):=\o X-L^{-1}\ell(\o X)\o \eta, \  or \ \ \phi:=I-L^{-1}
   \ell \otimes \o \eta,\vspace{-0.2cm}
\end{equation}
where $\ell$ is the  $\pi$-form  given  by
$\ell(\overline{X})=L^{-1}g(\overline{X},\overline{\eta}) ,$ then we
have{\em\,\!{:}}
\begin{description}
    \item[(a)] $\hbar(\o X,\o Y)=g(\phi(\o X), \o Y),$ \ \ \ \ \ \ \ \;\quad \textbf{{\em(}b{\em)}} $\phi(\o \eta)=0,$
    \ \ \ \ \ \ \ \ \ \quad \textbf{{\em(}c{\em)}} $\phi\,o\,\phi=\phi $,

   \item[(d)] $Tr(\phi)=n-1,$ \,\ \ \ \ \;\qquad\qquad\quad
    \textbf{{\em(}e{\em)}} $\nabla_{\beta \o X}\,\phi=0,$
     \ \ \ \ \ \ \ \quad \textbf{{\em(}f{\em)}}$\ \nabla_{\beta \o
     X}\,\hbar=0$.
\end{description}
\end{lem}

 As we have seen, a Landsberg
manifold is general Landsberg. The converse is not true.
Nevertheless, we have\vspace{-0.1cm}
\begin{prop}\label{pp.1a} A $C$-reducible general Landsberg
 manifold $(M,L)$ is a Landsberg manifold.
\end{prop}

\prof  Since $(M,L)$ is a $C$-reducible manifold, then, by
Definition \ref{3.def.1}, Lemma \ref{cor.i1}, the symmetry of
$\hbar$ and the  non-degeneracy of $g$, we get \vspace{-0.2cm}
\begin{equation*}\label{3.eq.1.1}
     T(\o X,\o Y) = \frac{1}{n+1}\{\hbar(\o X,\o Y) \o C
  +C(\o X) \phi(\o Y) +C(\o Y) \phi(\o X)\},\vspace{-0.3cm}
\end{equation*}
where   $\o C$ is the $\pi$-vector field defined by $g(\o C,\o
X):=C(\o X)$. Taking the $h$-covariant derivative $\nabla_{\beta \o
Z}$ of both sides of the above equation, we obtain\vspace{-0.2cm}
\begin{eqnarray*}\label{3.eq.1}
  (\nabla_{\beta \o Z}\,T)(\o X,\o Y) &=& \frac{1}{n+1}\{(\nabla_{\beta \o Z}\,\hbar)(\o X,\o Y) \o C+
  \hbar(\o X,\o Y)\nabla_{\beta \o Z}\,\o C+C(\o X) (\nabla_{\beta \o Z}\,\phi)(\o Y)+ \\
 &&+(\nabla_{\beta \o Z}\,C)(\o X)\phi(\o Y)+C(\o
 Y)(\nabla_{\beta \o Z}\,\phi)(\o X)+(\nabla_{\beta \o Z}\,C)(\o Y)\phi(\o X)\},\vspace{-0.2cm}
\end{eqnarray*}
 from which, by
setting $\o Z=\o \eta$ and taking into account the fact that
$\nabla_{\beta \o Z}\,\hbar=0$ and that $\nabla_{\beta \o
Z}\,\phi=0$\, ( Lemma \ref{cor.i1}), we get\vspace{-0.2cm}
\begin{equation*}\label{3.eq.2}
 (\nabla_{\beta \o \eta}\,T)(\o X,\o Y) = \frac{1}{n+1}\{
   \hbar(\o X,\o Y)\nabla_{\beta \o \eta}\,\o C+
  (\nabla_{\beta \o \eta}\,C)(\o X) \phi(\o Y)+(\nabla_{\beta \o \eta}\,C)(\o
   Y)\phi(\o X)\}.\vspace{-0.2cm}
\end{equation*}
Now, under the given assumption that the $(M,L)$ is a general
Landsberg manifold, then $\nabla_{\beta \o \eta}\,C=0$ (Definition
\ref{def.2a}) and hence $\nabla_{\beta \o \eta}\,\o C=0$. Hence
$\nabla_{\beta \o \eta}\,T=0$  and the result follows. \ \ $\Box$

\vspace{8pt}
 \par
Also, a Berwald manifold is  Landsberg. The converse is by no means
true,\linebreak although we have no counter-examples.  Finding  a
Landsberg manifold which is not Berwald is still an open problem.
Nevertheless, we have\vspace{-0.1cm}

\begin{prop}\label{1.pp.1a}{\em\cite{r48}} A $C$-reducible  Landsberg
 manifold $(M,L)$ is a Berwald \linebreak manifold.
\end{prop}

 Combining the above two Propositions, we obtain the more powerful result\,:\vspace{-0.1cm}
\begin{prop}\label{2.pp.1a} A $C$-reducible general  Landsberg
 manifold $(M,L)$ is a Berwald \linebreak manifold.
\end{prop}

 Summing up, we get\,\!:\vspace{-0.1cm}
\begin{thm}\label{3.th.1}Let $(M,L)$ be a $C$-reducible Finsler
manifold. The following assertion are
equivalent\,{\em:}\vspace{-0.2cm}
\begin{description}
    \item[(a)] $(M,L)$  is a Berwald manifold.
    \item[(b)] $(M,L)$  is a  Landsberg manifold.
    \item[(c)] $(M,L)$  is a general Landsberg manifold.
\end{description}
\end{thm}

We retrieve here a result of Matsumuoto {\cite{r2}},
namely\vspace{-0.2cm}
\begin{cor}\label{pp.2a}If the $h$-curvature tensor $R$  and $hv$-curvature tensor $P$ of a
 $C$-reducible manifold  vanish, then the manifold is Locally Minkowskian.
\end{cor}

\begin{rem}{\em\cite{r2}} It may be conjectured that a Finsler manifold will be
 Minkowskian if the $h$-curvature tensor $R$ and $hv$-curvature tensor $P$
 vanish. As above seen the conjecture is verified already under
 somewhat strong condition \lq\lq\,$C$-reducibility".
\end{rem}

 \begin{thm} \label{3.th.2}Let $(M,L)$ be a Finsler manifold. Then we have{\,\em:}
\begin{description}
\item[(a)] A $C$-reducible manifold is $P$-reducible.

   \item[(b)]  A $P$-reducible general Landsberg manifold  is  Landsberg.
\end{description}
\end{thm}

\prof\\
\textbf{(a)}  Since $(M,L)$ is $C$-reducible, then by Definition
\ref{3.def.1}, we have\vspace{-0.2cm}
\begin{equation*}
T(\o X,\o Y, \o Z)=\frac{1}{n+1}\mathfrak{S}_{\o X,\o Y,\o Z}\{
\hbar(\o X,\o Y)C(\o Z)\}.\vspace{-0.2cm}
\end{equation*}
Applying the $h$-covariant derivative $\nabla_{\beta \o W }$\  on
both sides of the above equation, taking into account the fact that
$(\nabla_{\beta \o W}\,T)(\o X,\o Y, \o Z)=g((\nabla_{\beta \o
W}\,T)(\o X,\o Y), \o Z)$ and that $\nabla_{\beta \o W}\,\hbar=0$,
we obtain
\begin{equation*}\label{3.eq.6}
g((\nabla_{\beta \o W}T)(\o X,\o Y), \o
Z)=\frac{1}{n+1}\mathfrak{S}_{\o X,\o Y,\o Z}\{\hbar(\o X,\o
Y)(\nabla_{\beta \o W}\,C)(\o Z)\}.\vspace{-0.2cm}
\end{equation*}
From which, by setting $\o W=\o \eta$ and noting that $P(\o X,\o
Y)\o \eta=(\nabla_{\beta \o \eta}\,T)(\o X,\o Y)$, the result
follows.\\
\textbf{(b)} Since $(M,L)$ is a $P$-reducible manifold, then
 by Definition \ref{def.p1}, taking into account the fact
  that $g$ is nondegenerate, we obtain  \vspace{-0.1cm}
\begin{equation}\label{3.eq.3}
P(\o X,\o Y)\o \eta=\delta(\o X)\phi(\o Y)+\delta(\o Y)\phi(\o X)+\o
\hbar(\o X,\o Y)\,\o \zeta ,\vspace{-0.1cm}
\end{equation}
where $\o \zeta$ is the $\pi$-vector field defined by $g(\o \zeta,
\o
X):=\delta(\o X)$.\\
 Since $\delta(\o \eta)=0$, then $Tr\{\o Y
\longmapsto \delta(\o Y)\phi(\o X)+\hbar(\o X,\o Y)\,\o \zeta \}=2
\delta(\o X)$. Taking the trace of both sides of (\ref{3.eq.3}),
using the fact that $P(\o X,\o Y)\o \eta=(\nabla_{\beta \o
\eta}\,T)(\o X,\o Y)$ (Lemma \ref{le.p}) and that $Tr\{\o Y
\longmapsto(\nabla_{\beta \o \eta}\,T)(\o X,\o Y)\}=(\nabla_{\beta
\o \eta}\,C)(\o X)$, we get\vspace{-0.2cm}
\begin{equation}\label{3.eq.4}
\delta(\o X)=\frac{1}{n+1}(\nabla_{\beta \o \eta}\,C)(\o
X).\vspace{-0.2cm}
\end{equation}
Now, from Equations (\ref{3.eq.3}) and (\ref{3.eq.4}), we
have\vspace{-0.2cm}
\begin{equation}\label{3.eq.5}
g(P(\o X,\o Y)\o \eta, \o Z)=\frac{1}{n+1}\mathfrak{S}_{\o X,\o Y,\o
Z}\{   \hbar(\o X,\o Y)(\nabla_{\beta \o \eta}\,C)(\o
Z)\}.\vspace{-0.2cm}
\end{equation}
According to the given assumption that the manifold is  general
Landsberg, then $\nabla_{\beta \o \eta}\,C=0$. Therefore, from
(\ref{3.eq.5}), we get  $P(\o X,\o Y)\o \eta=0$ and hence the
manifold is Landsberg.
 \ \ $\Box$
 \begin{prop}~\par\vspace{-0.2cm}
\begin{description}
    \item[(a)]A $C^{h}$-recurrent manifold is a $P^*$-Finsler manifold.
    \item[(b)]A general Landsberg $P^*$-Finsler manifold
    is a Landsberg manifold.
   \end{description}
\end{prop}
\prof The proof is straightforward and we omit it. \ \ $\Box$

\begin{prop}\label{733} A $C_{2}$-like Finsler manifold is a Berwald manifold
if, and only if, the $\pi$-tensor field $C$ is horizontally
parallel.
\end{prop}
\prof Let $(M,L)$ be  $C_{2}$-like. Then,
       $T(\o X,\o Y,\o Z)=\frac{1}{C(\o C)}C(\o X) C(\o Y) C(\o Z)$,
       from which
       $ T(\o X,\o Y)=\frac{1}{C(\o C)}C(\o X) C(\o Y) \o C.$
Taking the $h$-covariant derivative of both sides, we get
\begin{eqnarray*}
  (\nabla_{\beta \o Z}T)(\o X,\o Y)&=&\frac{-\nabla_{\beta \o Z}C(\o
C)}{C^{4}}C(\o X) C(\o Y) \o C+\frac{1}{C(\o C)}(\nabla_{\beta \o
Z}C)(\o X) C(\o Y) \o C+ \\
 &&+\frac{1}{C(\o C)}(\nabla_{\beta \o
Z}C)(\o Y) C(\o X) \o C+\frac{1}{C(\o C)}C(\o X) C(\o Y)
\nabla_{\beta \o Z}\o C.\vspace{-0.1cm}
\end{eqnarray*}
In view of this relation, $\nabla_{\beta \o Z}\,T=0$ if, and only
if, $\nabla_{\beta \o Z}\,C=0$. Hence the result. \ \ $\Box$

\vspace{-0.2cm}
\begin{cor} A $C_{2}$-like general Landsberg
 manifold  is a Landsberg manifold.\vspace{-1pt}
\end{cor}
 In view of the above Theorems, we
have\,\!:  \vspace{-0.2cm}
\begin{cor}\label{co.1}
The two notions of being Landsberg and general Landsberg coincide in
the case of $C$-reducibility, $P$-reducibility, $C_{2}$-likeness or
$P^*$-Finsler.
\end{cor}
As we know, a $C$-reducible Landsberg manifold is a Berwald manifold
(Proposition \ref{1.pp.1a} ). Moreover, A $C_{2}$-like Finsler
manifold is a Berwald manifold if, and only if, the $\pi$-tensor
field $C$ is horizontally parallel (Proposition \ref{733}).  We
shall try to generalize these results to the case of
semi-$C$-reduciblity.

\begin{thm} A semi-$C$-reducible Finsler manifold is a Berwald manifold
if, and only if, the characteristic scalar $\mu$ and the
$\pi$-tensor field $C$ are horizontally parallel.
\end{thm}

\prof Firstly, if  $(M,L)$ is  semi-$C$-reducible, then
\vspace{-0.1cm}
       $$T(\o X,\o Y,\o Z)=\frac{\mu}{n+1}\mathfrak{{S}}_{\o X,\o Y,\o Z}\{\hbar(\o X,\o Y)C(\o Z)\}
       +\frac{\tau}{C(\o C)}C(\o X) C(\o Y) C(\o Z).\vspace{-0.1cm}$$
Taking the $h$-covariant derivative of both sides, noting that
$\nabla_{\beta \o X}\hbar=0$, we get\vspace{-0.3cm}
\begin{eqnarray*}
  (\nabla_{\beta \o W}T)(\o X,\o Y,\o Z)&=&\frac{1}{n+1}
  \mathfrak{{S}}_{\o X,\o Y,\o Z}\{\hbar(\o X,\o Y)\{\mu(\nabla_{\beta \o
W}C)(\o Z)+(\nabla_{\beta \o W}\mu)C(\o Z)\}\}+\\
&& + \frac{\tau}{C^{2}}\mathfrak{{S}}_{\o X,\o Y,\o
Z}\{(\nabla_{\beta \o W}C)(\o X)C(\o Y)C(\o Z)\}-
\\
&& -\{\frac{\nabla_{\beta \o W}\,\mu}{C^{2}}+ \frac{\tau\,
\nabla_{\beta \o W}C(\o C)}{C^{4}}\}C(\o X)C(\o Y)C(\o Z)
 .\vspace{-0.1cm}
\end{eqnarray*}
\par Now, if the characteristic scalar $\mu$ and the  $\pi$-tensor
field $C$ are horizontally parallel, then $\nabla_{\beta \o W}T=0$
and $(M,L)$ is a Berwald manifold.
\par
Conversely, if $(M,L)$ is a Berwald manifold, then $ \nabla_{\beta
\o X}T=0$ and hence  $ \nabla_{\beta \o X}C=0$,  $ \nabla_{\beta \o
X}\o C=0$. These, together with the above equation, give
\begin{equation*}
\nabla_{\beta \o W}\mu \{\frac{1}{n+1}\mathfrak{{S}}_{\o X,\o Y,\o
Z}\{ \hbar(\o X,\o Y)C(\o Z)\}-\frac{1}{C^{2}}C(\o X)C(\o Y)C(\o
Z)\}=0,
\end{equation*}
 which implies immediately that $\nabla_{\beta \o W}\mu =0$.    \ \ $\Box$

\vspace{8pt}
 The following lemmas are useful for subsequent use\vspace{-0.2cm}
\begin{lem}\label{bracket}For all $\,\overline{X},\overline{Y}\in \mathfrak{X}(\pi(M))$,
 we have\,\em:\vspace{-0.2cm}
   \begin{description}
     \item[(a)] $[\gamma \overline{X},\gamma \overline{Y}]=
     \gamma(\nabla_{\gamma \overline{X}}\overline{Y}-
     \nabla_{\gamma \overline{Y}}\overline{X})$

     \item[(b)] $[\gamma \overline{X},\beta \overline{Y}]=-
     \gamma(P(\overline{Y},\overline{X})\overline{\eta}+\nabla_
     {\beta \overline{Y}}\overline{X})
     +\beta( \nabla_{\gamma \overline{X}}\overline{Y}-T(\overline{X},\overline{Y}))$

     \item[(c)] $[\beta \overline{X},\beta \overline{Y}]=
     \gamma(R(\overline{X},\overline{Y})\overline{\eta})
     + \beta(\nabla_{\beta \overline{X}}\overline{Y}-
     \nabla_{\beta \overline{Y}}\overline{X})$
 \end{description}
\end{lem}

\begin{lem}\label{lem}For all $\,\o {X},\o {Y}, \o Z, \o W \in \cp$ and $W\in\cpp$,
 we have\,\em:\vspace{-0.2cm}
\begin{description}
    \item[(a)]$g((\nabla_{W}T)(\o X,\o Y),\o
Z)=g((\nabla_{W}T)(\o X,\o Z),\o Y)$,
    \item[(b)]$ g(S(\o X,\o Y)\o Z, \o
W)=-g(S(\o X,\o Y)\o W,\o Z)$.
\end{description}
\end{lem}

\prof\\
\textbf{(a)} From the definition of the covariant derivative, we
get\vspace{-0.2cm}
\begin{equation}\label{11}
\left.
    \begin{array}{rcl}
g((\nabla_{W}T)(\o X,\o Y),\o Z)&=&g(\nabla_{W}T(\o X,\o Y),\o
Z)-g(T(\nabla_{W}\o X,\o Y),\o Z)-\\ &&-g(T(\o X,\nabla_{W}\o Y),\o
Z). \vspace{-0.2cm}
 \end{array}
  \right.
\end{equation}
Now, we have\vspace{-0.2cm}
\begin{equation*}
\left.
    \begin{array}{rcl}
  g(\nabla_{W}T(\o X,\o Y),\o Z)&=&W \cdot g(T(\o X,\o Y),\o Z)- g(T(\o X,\o Y),\nabla_{W}\o Z)\\
   &=&W \cdot g(T(\o X,\o Y),\o Z)- g(T(\o X,\nabla_{W}\o Z),\o Y),\vspace{-0.2cm}
 \end{array}
  \right.
\end{equation*}
Similarly,\vspace{-0.2cm}
\begin{equation*}
\left.
    \begin{array}{rcl}
  g(T(\o X,\nabla_{W}\o Y),\o Z)&=&
  W \cdot g(T(\o X,\o Z),\o Y)- g(\nabla_{W}T(\o X,\o Z),\o Y).\vspace{-0.2cm}
 \end{array}
  \right.
\end{equation*}
 Substituting  these two equations into (\ref{11}), noting the property that
$g(T(\nabla_{W}\o X,\o Y),\o Z)$\\ $= g(T(\nabla_{W}\o X,\o Z),\o
Y)$ (cf. \S 1), the result follows.

\vspace{4pt} \noindent
\textbf{(b)} follows directly from the
general formula (which can be easily proved)\vspace{-0.2cm}
\begin{equation*}
 g(\textbf{K}( X, Y)\o Z, \o
W)+g(\textbf{K}( X, Y)\o W,\o Z)=0\vspace{-0.2cm}
\end{equation*}
 by setting
$X=\gamma \o X$ and $Y=\gamma \o Y$, where $\textbf{K}$ is the
classical curvature tensor of the Cartan connection  as a linear
connection in the pull-back bundle (cf. \S 1).  \ \ $\Box$

\begin{prop}\label{7p}Let $(M,L)$ be  a $C^h$-recurrent Finsler manifold
{\em{(}$\nabla_{\beta\o X}T=\lambda_{0}(\o X)T$)}. Then, we
have{\,\!\em{:}}
\begin{description}
    \item[(a)] If $K_{o}:=\lambda_{o}(\o \eta)=0$, then the
    $hv$-curvature tensor $P$ is expressed in the form\,\!{\em :}\\
      $ P(\o X,\o Y,\o Z, \o W)=\lambda_{o}(\o Z)T(\o X,\o Y,\o W)-\lambda_{o}(\o W)T(\o X,\o Y,\o
    Z)$\\
    and the $(v)hv$-torsion $\widehat{P}$ vanishes.

    \item[(b)] If $K_{o}\neq0$, then   the $v(hv)$-torsion tensor $\widehat{P}$ is recurrent\,\!{\em :}\\
    $( \nabla_{\beta \o Z} \widehat{P})(\o X,\o Y)=(\lambda_{o}(\o Z)+
    \frac{\nabla_{\beta \o Z} K_{o}}{K_{o}})\widehat{P}(\o X,\o Y)$.
\end{description}
\end{prop}

\prof\\
\textbf{(a)} The $hv$-curvature tensor $P$ can be written in the
form \cite{r48}:\vspace{-0.1cm}
\begin{equation*}\label{7a6}
    \begin{array}{rcl}
   P(\o X,\o Y,\o Z, \o W) &=&g(( \nabla_{\beta \o Z}T)(\o X, \o Y),\o W)
   -g(( \nabla_{\beta \o W}T)(\o X, \o Y),\o Z)+\\
   &&+ g(T(\o X,\o Z),\widehat{P}(\o W,\o Y))-g(T(\o X,\o W),\widehat{P}(\o Z,\o
   Y)).\vspace{-0.1cm}
\end{array}
\end{equation*}
Then, by using $\widehat{P}(\o X,\o Y)=( \nabla_{\beta \o \eta}T)(\o
X, \o Y)$ (Lemma \ref{le.p}) and the $C^h$-recurrence condition, we
get\vspace{-0.1cm}
\begin{equation*}\label{7a4}
  \left.
    \begin{array}{rcl}
   P(\o X,\o Y,\o Z, \o W) &=&\lambda_{o}(\o Z)T(\o X,\o Y,\o W)-\lambda_{o}(\o W)T(\o X,\o Y,\o
    Z)-\\
    &&-\lambda_{o}(\o \eta)\{g(T(\o X,\o W),T(\o Y,\o Z))
    -g(T(\o X,\o Z),T(\o Y,\o W))\}\\
    &=&\lambda_{o}(\o Z)T(\o X,\o Y,\o W)-\lambda_{o}(\o W)T(\o X,\o Y,\o
    Z)-\lambda_{o}(\o \eta)S(\o X,\o Y,\o Z,\o W).\vspace{-0.1cm}
\end{array}
  \right.
\end{equation*}
Now, if $\lambda_{o}(\o \eta)=0$, then (a) follows from the above
relation.

\vspace{4pt} \noindent \textbf{(b)} If $K_{o}:=\lambda_{o}(\o
\eta)\neq0$, then by Lemma \ref{le.p} and the recurrence condition,
we have\vspace{-0.1cm}
\begin{equation*}\label{72}
\widehat{P}(\o X,\o Y)=K_{o}T(\o X, \o Y),
\end{equation*}
from which
  $$( \nabla_{\beta \o Z}\widehat{P})(\o X,\o Y)=\{\nabla_{\beta \o Z}K_{o}+K_{o}\lambda_{o}(\o Z)\}T(\o X,\o
  Y).$$
Then, (b) follows from the above two equations.
    \ \  $\Box$

\begin{thm}\label{th.a1}Assume that $(M,L)$ is $C^{h}$-recurrent. Then, the
$v$-curvature tensor $S$ is recurrent with respect to the
$h$-covariant differentiation\;{\em{:}} $\nabla_{\beta \o
X}S=\theta(\o X)S$, where $\theta$ is a $\pi$-form  of order one.
\end{thm}

\prof  One can  easily show that\,: For all $X,Y,Z \in
\cpp$,\vspace{-0.2cm}
\begin{equation*}\label{3.eq.10}
    \mathfrak{S}_{X,Y,Z}\{  \textbf{K}(X,Y)\rho
    Z+\nabla_{X}\textbf{T}(Y,Z)+\textbf{T}(X,[Y,Z])\}=0.\vspace{-0.2cm}
\end{equation*}
Setting $X=\gamma \o X,\  Y=\gamma \o Y$ and $Z=\beta \o Z$ in the
above equation, we get\vspace{-0.1cm}
\begin{eqnarray*}
   S(\o X, \o Y)\o Z&=&\nabla_{\gamma \o Y}T(\o X,\o Z )-\nabla_{\gamma \o X}T(\o Y,\o Z )
   -\nabla_{\beta \o Z}\textbf{T}(\gamma \o X,\gamma \o Z ) -\\
   &&-\textbf{T}(\gamma \o X,[\gamma \o Y,\beta \o Z]) +\textbf{T}(\gamma \o Y,[\gamma \o X,\beta \o Z])
   +\textbf{T}([\gamma \o X,\gamma \o Y],\beta \o Z). \vspace{-0.1cm}
\end{eqnarray*}
 Using Lemma \ref{bracket} and the fact
 that $\textbf{T}(\gamma \o X,\gamma \o Z )=0$, the above equation  reduces to
\begin{equation}\label{La}
\left.
    \begin{array}{rcl}
S(\o X,\o Y)\o Z &=& (\nabla_{\gamma \o Y}T)(\o X,\o
Z)-(\nabla_{\gamma \o X}T)(\o Y,\o Z)+\\
& &+T(\o X,T(\o Y,\o Z))-T(\o Y,T(\o X,\o Z)).\vspace{-0.2cm}
    \end{array}
  \right.
\end{equation}
From which, since $g(T(\o X,\o Y),\o Z)=g(T(\o X,\o Z),\o Y)$, we
have\vspace{-0.2cm}
\begin{equation*}\label{3.eq.7}
\left.
    \begin{array}{rcl}
g(S(\o X,\o Y)\o Z,\o W) &=& g((\nabla_{\gamma \o Y}T)(\o X,\o
Z), \o W)-g((\nabla_{\gamma \o X}T)(\o Y,\o Z),\o W)+\\
& &+g(T(\o X,\o W),T(\o Y,\o Z))-g(T(\o Y,\o W),T(\o X,\o
Z)).\vspace{-0.2cm}
    \end{array}
  \right.
\end{equation*}
Similarly,
\begin{equation*}\label{3.eq.7}
\left.
    \begin{array}{rcl}
g(S(\o X,\o Y)\o W,\o Z) &=& g((\nabla_{\gamma \o Y}T)(\o X,\o
W), \o Z)-g((\nabla_{\gamma \o X}T)(\o Y,\o W),\o Z)+\\
& &+g(T(\o X,\o Z),T(\o Y,\o W))-g(T(\o Y,\o Z),T(\o X,\o
W)).\vspace{-0.2cm}
    \end{array}
  \right.
\end{equation*}
The above two equations, together with  Lemma \ref{lem},
yield\vspace{-0.1cm}
\begin{equation}\label{3.eq.8}
    g((\nabla_{\gamma \o X}T)(\o Y,\o Z),\o
W)=g((\nabla_{\gamma \o Y}T)(\o X,\o Z),\o W).\vspace{-0.1cm}
\end{equation}
By (\ref{La}) and (\ref{3.eq.8}), we obtain \vspace{-0.2cm}
 \begin{equation}\label{3.eq.9}
S(\o X, \o Y , \o Z, \o W)=g(T(\o X, \o W) , T(\o Y, \o Z))-g(T(\o
Y, \o W) ,T( \o X, \o Z)).\vspace{-0.2cm}
\end{equation}
\par
Now, using the given assumption that the manifold is
$C^{h}$-recurrent, Equation (\ref{3.eq.9}) implies
that\vspace{-0.2cm}
\begin{eqnarray*}
   (\nabla_{\beta \o X}S)(\o Y , \o Z, \o V, \o W )&=&\nabla_{\beta \o X}S(\o Y , \o Z, \o V, \o W
   )-\\
   &&-S(\nabla_{\beta \o X}\o Y , \o Z, \o V, \o W )-S(\o Y , \nabla_{\beta \o X}\o Z, \o V, \o W )-\\
   &&-S(\o Y , \o Z, \nabla_{\beta \o X}\o V, \o W )-S(\o Y , \o Z, \o V, \nabla_{\beta \o X}\o W
   ).\\
   &=&+\nabla_{\beta \o X}g(T(\o Y, \o W) , T(\o Z, \o V))-
   \nabla_{\beta \o X}g(T(\o Z, \o W) , T(\o Y, \o V))-\\
    &&-g(T(\nabla_{\beta \o X}\o Y, \o W) , T(\o Z, \o V))+
   g(T(\o Z, \o W) , T(\nabla_{\beta \o X}\o Y, \o V))-\\
   &&-g(T(\o Y, \o W) , T(\nabla_{\beta \o X}\o Z, \o V))+
   g(T(\nabla_{\beta \o X}\o Z, \o W) , T(\o Y, \o V))-\\
    &&-g(T(\o Y, \o W) , T(\o Z, \nabla_{\beta \o X}\o V))+
   g(T(\o Z, \o W) , T(\o Y, \nabla_{\beta \o X}\o V))-\\
   &&-g(T(\o Y, \nabla_{\beta \o X}\o W) , T(\o Z, \o V))+
   g(T(\o Z, \nabla_{\beta \o X}\o W) , T(\o Y, \o V)).\\
    &=&g((\nabla_{\beta \o X}T)(\o Y, \o W) , T(\o Z, \o V))+
   g(T(\o Y, \o W) , (\nabla_{\beta \o X}T)(\o Z, \o V))-\\
   & &- g((\nabla_{\beta \o X}T)(\o Z, \o W) , T(\o Y, \o V))-
   g(T(\o Z, \o W) , (\nabla_{\beta \o X}T)(\o Y, \o V)).\\
   &=&2 \lambda_{o}(\o X) S(\o Y , \o Z, \o V, \o W )=:\theta(\o X)S(\o Y , \o Z, \o V, \o W
   ).\vspace{-0.2cm}
\end{eqnarray*}
Hence, the result follows. \ \ $\Box$

\begin{cor}\label{CC} In the course of the proof of Theorem {\em\ref{th.a1}}, we have
shown that {\em{(Equations (\ref{3.eq.8}) and
(\ref{3.eq.9}))}}{\em\,:}
\begin{description}
    \item[(a)] $ (\nabla_{\gamma \o X}T)(\o Y,\o Z)=(\nabla_{\gamma \o Y}T)(\o X,\o Z) $,

    \item[(b)]$S(\o X, \o Y , \o Z, \o W)=g(T(\o X, \o W) , T(\o Y, \o Z))-g(T(\o
Y, \o W) ,T( \o X, \o Z))$.
   \end{description}
\end{cor}

\begin{cor}Let $(M,L)$ be a $C_{2}$-like Finsler manifold. Then the
the $v$-curvature tensor $S$ vanishes.
\end{cor}
\prof Substituting \,$ T(\o X,\o Y)=\frac{1}{C(\o C)}C(\o X) C(\o Y)
\o C$\, in Corollary \ref{CC}(b), we get the result. \ \ $\Box$

\begin{cor}Let $(M,L)$ be a $C$-reducible manifold. Then,
\begin{description}
    \item[(a)] the $v$-curvature tensor $S$ has the form\vspace{-0.2cm}
   \begin{eqnarray*}
     S(\o X,\o Y,\o Z,\o W) &=& \frac{1}{(n+1)^2}\{C^{2}\hbar (\o X,\o W)
     \hbar(\o Y,\o Z)-C^{2}\hbar (\o Y,\o W)
     \hbar(\o X,\o Z)+\\
     &&+\hbar (\o X,\o W)C(\o Y)C(\o Z) +\hbar (\o Y,\o Z)C(\o X)C(\o W)-\\
     &&- \hbar (\o Y,\o W)C(\o X)C(\o Z)
      -\hbar (\o X,\o Z)C(\o Y)C(\o W)\}.
   \end{eqnarray*}
    \item[(b)]the vertical Ricc  tensor $Ric^{v}$ has the
    form\vspace{-0.2cm}

    $$ Ric^{v}(\o X,\o Y)=\frac{(3-n)}{(n+1)^{2}}C(\o X)C(\o
    Y)-\frac{(n-1)}{(n+1)^{2}}\,
    C^{2}\hbar(\o X,\o Y).$$
    \item[(c)]the vertical scalar curvature $Sc^{v}$ has the
    form\vspace{-0.2cm}
     $$ Sc^{v}=\frac{(2-n)}{(n+1)}\,C^{2}.$$
\end{description}
\end{cor}

\begin{thm}\label{th.a2} A Finsler manifold $(M,L)$ is $P$-Symmetric if, and only if,
the\linebreak $v$-curvature tensor $S$ satisfies the equation
$\nabla_{\beta \o \eta}S=0$.
\end{thm}

\prof One can show that: For all $X,Y,Z \in \cpp$,\vspace{-0.1cm}
\begin{equation}\label{K}
    \mathfrak{S}_{X,Y,Z}\{  \nabla_{Z}\textbf{K}(X,Y)-\textbf{K}(X,Y)\nabla_{Z}
    -\textbf{K}([X,Y],Z)\}=0.\vspace{-0.1cm}
\end{equation}
Setting $X=\gamma \o X, Y=\gamma \o Y$ and $Z=\beta \o Z$ in the
above equation, we get\vspace{-0.1cm}
\begin{equation*}
\left.
    \begin{array}{rcl}
&&\nabla_{\beta\o Z}S(\o X,\o Y)\o W + \nabla_{\gamma \o Y}P(\o
Z,\o X) \o W-\nabla_{\gamma \o X}P(\o Z,\o Y)\o W-\\
& &- S(\o X,\o Y)\nabla_{\beta\o Z}\o W+P(\o Z,\o Y) \nabla_{\gamma
\o X}\o W-P(\o Z,\o X) \nabla_{\gamma
\o Y}\o W-\\
& &-\textbf{K}([\gamma\o X ,\gamma\o Y],\beta \o Z)\o
W-\textbf{K}([\gamma\o Y ,\beta\o Z],\gamma\o X)\o
W-\textbf{K}([\beta\o Z ,\gamma\o X],\gamma \o Y)\o W=0.
\vspace{-0.1cm}
    \end{array}
  \right.
\end{equation*}
By using  Lemma \ref{bracket}, the above relation reduces to
\vspace{-0.1cm}
\begin{equation}\label{3.eq.11}
\left.
    \begin{array}{rcl}
&&(\nabla_{\beta\o Z}S)(\o X,\o Y,\o W )+ (\nabla_{\gamma \o Y}P)(\o
Z,\o X, \o W)-(\nabla_{\gamma \o X}P)(\o Z,\o Y,\o W)+\\& &+S(P(\o
Z,\o Y)\o \eta,\o X)\o W- S(P(\o Z,\o X)\o \eta,\o Y)\o W +\\& &+
P(T(\o Y,\o Z),\o X)\o W-P(T(\o X,\o Z),\o Y)\o W=0. \vspace{-0.1cm}
    \end{array}
  \right.\vspace{-0.2cm}
\end{equation}
Setting $\o Z=\o \eta$\; in the above equation, taking
 into account Lemma \ref{le.p} and the fact that
$T(\o X,\o \eta)=0$ and that $(\nabla_{\gamma \o X}P)(\o \eta,\o
Y,\o Z)=-P(\o X,\o Y)\o Z$, we get\vspace{-0.1cm}
\begin{equation}\label{3.eq.12}
    P(\o X, \o Y ) \o Z=P(\o Y, \o X ) \o Z -(\nabla_{\beta \o \eta}S)(\o X, \o Y, \o Z).\vspace{-0.1cm}
\end{equation}
The result follows immediately from (\ref{3.eq.12}). \ \ $\Box$

\vspace{8pt}
 According to (\ref{3.eq.12}) and Lemma \ref{le.p}, we have\,:\vspace{-0.2cm}
\begin{cor}\label{co.2}Let $\widehat{P}(\o X,\o Y):=P(\o X,\o Y)\o \eta$ and $\widehat{T}(\o X,\o
Y):=(\nabla_{\beta\o \eta}T)(\o X,\o Y)$.  Then the $\pi$-tensor
fields $\widehat{P}$ and $\widehat{T}$ are symmetric.
\end{cor}

 Theorem \ref{th.a1} and Theorem \ref{th.a2} give rise the following
 result.
\vspace{-0.2cm}
\begin{thm}\label{th.4}Assume that a Finsler manifold $(M,L)$ is $C^{h}$-recurrent and $P$-symmetric.
If $\theta(\o \eta)\neq0$, then the $v$-curvature tensor $S$
vanishes identically.
\end{thm}

Now, we shall prove the following lemma which provides some
important and useful properties of the torsion tensor $T$ and the
$v$-curvature $S$\,:\vspace{-0.2cm}
\begin{lem}\label{lem.1}For every $\o X, \o Y , \o Z \ and \ \o W\in \cp$, we have
\begin{description}

 \item[(a)] $T(\o X, \o Y )=T(\o Y, \o  X)$,
  \item[(b)] $T(\o \eta, \o  X)=0$,
 \item[(c)] $\mathfrak{S}_{\o X,\o Y,\o Z}S(\o X, \o Y )\o Z=0$,
    \item[(d)] $g(S(\o X, \o Y )\o Z, \o W)=g(S(\o Z, \o W )\o X, \o
    Y)$,
    \item[(e)]$S(\o \eta, \o X )\o Y=0=S(\o X, \o \eta )\o Y$,
    \item[(f)]$(\nabla_{\gamma\o X }S)(\o \eta, \o Y )\o Z=-S(\o X, \o Y )\o
    Z$, \  $(\nabla_{\gamma\o X }S)(\o \eta, \o X )\o \eta=0$ .
\item[(g)]$S(\o X,\o Y)\o
Z=-\frac{1}{2}\{(D_{\gamma
\overline{X}}T)(\overline{Y},\overline{Z})-(D_{\gamma
\overline{Y}}T)(\overline{X},\overline{Z})\}$.\\
Consequently, $S$ vanishes if and only if $\, (D_{\gamma
\overline{X}}T)(\overline{Y},\overline{Z})=(D_{\gamma
\overline{Y}}T)(\overline{X},\overline{Z})$. 
\end{description}
\end{lem}

\prof\\
\noindent\textbf{(a)} From Corollary \ref{CC}(a), we
have\vspace{-0.2cm}
 $$ (\nabla_{\gamma \o X}T)(\o Y,\o Z)=(\nabla_{\gamma \o Y}T)(\o X,\o Z).\vspace{-0.2cm}$$
 Setting $\o Z=\o \eta$ and using  the fact that $T(\o X,\o
 \eta)=0$ and that $K\,o\,\gamma =id_{\cp}$, the result
 follows.

 \vspace{5pt}
\noindent \textbf{(b)} Follows from  (a) together with the relation
$T(\o X,\o
 \eta)=0$.

 \vspace{5pt}
\noindent \textbf{(c)} Setting $ X=\gamma \o X,\ Y=\gamma \o Y $ and
$Z=\gamma \o Z$\; in
 (\ref{K}) and using Lemma \ref{bracket},  we get\vspace{-0.2cm}
 $$\mathfrak{S}_{\o X,\o Y,\o Z}(\nabla_{\gamma \o X}S)(\o Y, \o Z, \o W)=0.\vspace{-0.2cm}$$
Again, setting $\o W=\o \eta$ in the above equation and using the
fact that $S(\o X, \o Y )\o \eta=0$ and that $ K\,o\,
\gamma=id_{\cp}$, the result follows.

\vspace{5pt}
\noindent\textbf{(d)} Follows from Corollary
\ref{CC}(b), noting that $T$ is symmetric.

\vspace{5pt}
\noindent\textbf{(e)} and \textbf{(f)} are clear.

\vspace{5pt}
\noindent\textbf{(g)} From the relation \,$D_{\gamma
\overline{X}}\o Y=\nabla_{\gamma \overline{X}}\o Y-T(\o X,\o Y)$
\cite{r62}, we get
\begin{equation*}
      (D_{\gamma \overline{X}}T)(\o Y,\o Z)= (\nabla_{\gamma \overline{X}}T)(\o Y,\o Z)
   -T(\o X,T(\o Y,\o Z))
    +T(T(\o X, \o Y),\o Z)+T(\o Y,T(\o X,\o Z)),
\end{equation*}
\begin{equation*}
      (D_{\gamma \overline{Y}}T)(\o X,\o Z)= (\nabla_{\gamma \overline{Y}}T)(\o X,\o Z)
   -T(\o Y,T(\o X,\o Z))
    +T(T(\o Y, \o X),\o Z)+T(\o X,T(\o Y,\o Z)).
   \end{equation*}
The result follows from the above two equations, using Corollary
\ref{CC} and the symmetry of $T$.\ \ $\Box$
\vspace{0.2cm}
\par As a direct consequence of the above lemma, we
have the
\begin{cor}
\vspace{-0.2cm} A $P_{2}$-like Finsler manifold is $P$-symmetric.
\end{cor}

\begin{prop}\label{7th.a1}Assume that $(M,L)$ is $C^{v}$-recurrent. Then, the
$v$-curvature tensor $S$ is $v$-recurrent\;\!{\em{:}}
$\nabla_{\gamma\o X}S=\Psi(\o X)S$, $\Psi$ being a $(1)\pi$-form.
Consequently, $S$ vanishes identically.
\end{prop}
\prof Taking the $v$-covariant derivative of both sides of the
relation in Corollary \ref{CC}(b)\, and, then, using the assumption
that $\nabla_{\gamma \overline{X}}T=\lambda_{0}(\overline{X})T$, we
get \vspace{-0.1cm}
\begin{equation*}
  (\nabla_{\gamma \o X}S)(\o Y , \o Z, \o V, \o W )
       =2 \lambda_{o}(\o X) S(\o Y , \o Z, \o V, \o W )=:\psi(\o X)S(\o Y , \o Z, \o V, \o
       W),
\vspace{-0.1cm}
\end{equation*} which shows that $S$ is $v$-recurrent.
\par Now, setting $\o V=\o \eta$ in the last equation, using the properties of $S$ and noting that
$K\,o\,\gamma=id_{\cp}$, we conclude that $S=0$. \ \ $\Box$
\vspace{7pt}
\par The following result gives a
characterization of Riemannian manifolds in terms of
$C^v$-recurrence and $C^0$-recurrence.\vspace{-0.2cm}
\begin{thm}~\par\vspace{-0.2cm}
\begin{description}
    \item[(a)] A $C^v$-recurrent Finsler manifold is Riemannian,
    \item[(b)] A $C^0$-recurrent Finsler manifold is Riemannian.
\end{description}
\end{thm}

\prof
    (a) Since $(M,L)$ is $C^v$-recurrent,
then $(\nabla_{\gamma \overline{X}}T)(\overline{Y},\overline{Z})
 =\lambda_{o} (\overline{X}) T(\overline{Y},\overline{Z}),$ from
 which, by setting $\o X=\o \eta$ and noting that $\nabla_{\gamma\o\eta}T=-T$, we get\vspace{-0.1cm}
\begin{equation}\label{7c.1}
T(\overline{Y},\overline{Z})
 =-\lambda_{o} (\overline{\eta}) T(\overline{Y},\overline{Z}).\vspace{-0.1cm}
\end{equation}
But since $ (\nabla_{\gamma \o X}T)(\o Y,\o Z)=(\nabla_{\gamma \o
Y}T)(\o X,\o Z) $ (Corollary \ref{CC}), then  $ \lambda_{o}( \o
X)T(\o Y,\o
Z)=\lambda_{o}( \o Y)T(\o X,\o Z) $. Hence,\\
\vspace{-0.1cm}
\begin{equation}\label{7c.2}
\lambda_{o} (\overline{\eta})
T(\overline{Y},\overline{Z})=0.\vspace{-0.1cm}
\end{equation}
Then, the result follows from (\ref{7c.1}) and (\ref{7c.2}).\\
(b) can be proved similarly.\ \ $\Box$

\begin{thm}\label{th.5}For a Finsler manifold $(M,L)$, the following
assertions are \linebreak equivalent\,{\em:}\vspace{-0.2cm}
\begin{description}
    \item[(a)]  $(M,L)$ is $S^{v}$-recurrent.
    \item[(b)] The $v$-curvature tensor $S$ vanishes identically.
    \item[(c)]  $(M,L)$ is $S^{v}$-recurrent of the second order.
\end{description}
\end{thm}

\prof\\
$\textbf{(a)}\Longrightarrow\textbf{(b)}$\,: If $(M,L)$ is
$S^{v}$-recurrent, then  by Definition \ref{def.a1}(a) we have
\vspace{-0.2cm}
\begin{equation*}
   (\nabla_{\gamma\o W}S)(\o X, \o Y ,\o Z)=\lambda (\o W)S(\o Y, \o X ) \o Z,\vspace{-0.2cm}
\end{equation*}
from which,  by setting $\o Z=\o \eta$, taking into account the fact
that $S(\o X,\o Y)\o \eta=0$ and that $K o
\gamma=id_{\p}$, the result follows.\\
$\textbf{(b)}\Longrightarrow\textbf{(a)}$\,: Trivial.\\
$\textbf{(b)}\Longrightarrow\textbf{(c)}$\,: Trivial.\\
$\textbf{(c)}\Longrightarrow\textbf{(b)}$\,: If the given manifold
$(M,L)$ is $S^{v}$-recurrent of the second order, then by Definition
\ref{def.a1}(b)  we get \vspace{-0.2cm}
\begin{equation*}\label{eq.a4}
  \left.
    \begin{array}{rcl}
        \Theta(\o X,\o Y)S(\o Z, \o V)\o W& = &
        (\stackrel{2}\nabla\stackrel{2}\nabla S)(\o Y, \o X, \o Z, \o V,\o
        W)\\
        &=&  \nabla_{\gamma \o Y}(\nabla_{\gamma \o X}S)(\o Z, \o V,\o W)  -
  (\nabla_{\gamma \nabla_{\gamma \o Y}\o X }S)(\o Z, \o V,\o W) - \\
    &&-(\nabla_{\gamma \o X}S)(\nabla_{\gamma \o Y}\o Z, \o V,\o W)
   -(\nabla_{\gamma \o X}S)(\o Z, \nabla_{\gamma \o Y}\o V,\o W)- \\
   &&-(\nabla_{\gamma \o X}S)(\o Z, \o V,\nabla_{\gamma \o Y}\o W).
    \end{array}
  \right. \vspace{-0.2cm}
\end{equation*}
 By substituting $\o Z=\o \eta=\o W$ in the above equation  and using  Lemma \ref{lem.1}
 and  the fact that
$S(\o X,\o Y)\o \eta=0$, we get  $$S(\o X, \o Y)\o Z=-S(\o Z, \o
Y)\o X \ \  and  \ \   S(\o X, \o Y)\o Z=-S(\o X, \o Z)\o Y.$$ From
this, together with the identity $ \mathfrak{S}_{\o X,\o Y,\o Z}S(\o
X, \o Y)\o Z=0$, the $v$-curvature tensor $S$ vanishes identically.
\ \ $\Box$

\vspace{7pt} In view of the above theorem we have\,: \vspace{-0.2cm}
\begin{cor}\label{co.3}~\par \vspace{-0.1cm}
\begin{description}
    \item[(a)]An $S^{v}$-recurrent {\em({\it resp. a second order $S^{v}$-recurrent})}
    manifold $(M,L)$ is $S_{3}$-like, provided that $dim\,M\geq4$.
    \item[(b)]An $S^{v}$-recurrent {\em({\it resp. a second order $S^{v}$-recurrent})}
     manifold $(M,L)$ is $S_{4}$-like, provided that $dim\,M\geq5$.
\end{description}
\end{cor}

\begin{thm}If $(M,L)$ is a $P_{2}$-like Finsler manifold, then the
 $v$-curvature tensor $S$ vanishes or the $hv$-curvature tensor $P$ vanishes.
In the later case, the $h$-covariant derivative of $S$ vanishes.
\end{thm}

\prof As $(M,L)$ is $P_{2}$-like, then
$P(\overline{X},\overline{Y},\overline{\eta}, \o W)=
 \alpha(\overline{\eta})T (\overline{X},\overline{Y},\o W)=:
 \alpha_{o}T (\overline{X},\overline{Y},\o W)$ and hence
 \begin{equation}\label{7a5}
    \widehat{P}(\o X,\o Y)=\alpha_{o}T (\overline{X},\overline{Y}).
\end{equation}
Now, setting $\o W=\o \eta$ into (\ref{3.eq.11}), we
get\vspace{-0.1cm}
\begin{equation*}
\left.
    \begin{array}{rcl}
&& (\nabla_{\gamma \o Y}\widehat{P})(\o Z,\o X)-(\nabla_{\gamma \o
X}\widehat{P})(\o Z,\o Y)-P(\o Z,\o X)\o Y+P(\o Z,\o Y)\o X-\\
& &-\widehat{P}(T(\o X,\o Z),\o Y)+ \widehat{P}(T(\o Y,\o Z),\o
X)=0.
    \end{array}
  \right.\vspace{-0.1cm}
\end{equation*}
Hence, \vspace{-0.1cm}
\begin{equation*}
\left.
    \begin{array}{rcl}
&& g((\nabla_{\gamma \o Y}\widehat{P})(\o Z,\o X),\o
W)-g((\nabla_{\gamma \o X}\widehat{P})(\o Z,\o Y),\o W)-P(\o Z,\o
X,\o Y,\o W)+\\& &+P(\o Z,\o Y,\o X,\o W)-g(\widehat{P}(T(\o X,\o
Z),\o Y),\o W)+ g(\widehat{P}(T(\o Y,\o Z),\o X),\o W)=0.
    \end{array}
  \right.\vspace{-0.1cm}
\end{equation*}
From which, together with (\ref{7a5}) and Definition
 \ref{def.p2like}, taking into account the relation
 $ (\nabla_{\gamma \o Y}\widehat{P})(\o
Z,\o X)=(\nabla_{\gamma \o Y}\alpha_{o})T(\o Z,\o X)+
\alpha_{o}(\nabla_{\gamma \o Y}T)(\o Z,\o X)$, we obtain
\vspace{-0.1cm}
\begin{equation*}
\left.
    \begin{array}{rcl}
&& g((\nabla_{\gamma \o Y}\alpha_{o})T(\o Z,\o
X)+\alpha_{o}(\nabla_{\gamma \o Y}T)(\o Z,\o X),\o
W)-g((\nabla_{\gamma \o X}\alpha_{o})T(\o Z,\o
Y)+\\&&+\alpha_{o}(\nabla_{\gamma \o X}T)(\o Z,\o Y),\o W) +
\alpha(\overline{X})T (\overline{Z},\overline{Y},\o W)
 -\alpha(\overline{W})\, T(\overline{Z},\o Y,\overline{X})
 - \alpha(\overline{Y})T (\overline{Z},\overline{X},\o W)\\&&
 +\alpha(\overline{W})\, T(\overline{X},\o Y,\overline{Z})
-g(\alpha_{o}T(T(\o X,\o Z),\o Y),\o W)+ g(\alpha_{o}T(T(\o Y,\o
Z),\o X),\o W)=0.
   \end{array}
  \right.\vspace{-0.1cm}
\end{equation*}
Therefore, using Corollary \ref{CC},\vspace{-0.1cm}
\begin{equation*}
(\nabla_{\gamma \o Y}\alpha)(\o \eta)T(\o X,\o Z,\o W)-
(\nabla_{\gamma \o X}\alpha)(\o \eta)T(\o Y,\o Z,\o
W)=\alpha_{o}S(\o X,\o Y,\o W,\o Z).\vspace{-0.1cm}
\end{equation*}
It is to be observed that the left-hand side of the above equation
is symmetric in the arguments $\o Z$ and $\o W$ while the right-hand
side is skew-symmetric in the same arguments. Hence we
have\vspace{-0.1cm}
\begin{equation}\label{8a11}
\alpha_{o}S(\o X,\o Y,\o W,\o Z)=0,
\end{equation}
\vspace{-0.4cm}
\begin{equation}\label{7a11}
    \varepsilon(\o Y)T(\o X,\o Z,\o W)-
\varepsilon(\o X)T(\o Y,\o Z,\o W)=0,\vspace{-0.1cm}
\end{equation}
where $\varepsilon$ is the $\pi$-form defined by $\varepsilon(\o
Y):=(\nabla_{\gamma \o Y}\alpha)(\o \eta)$.
 \par Now, If
$\varepsilon\neq0$, it follows from (\ref{7a11}) that there exists a
scalar function $\Upsilon$ such that $T(\o X,\o Y,\o
Z)=\Upsilon\,\varepsilon(\o X) \varepsilon(\o Y) \varepsilon(\o Z)$.
Consequently,  $ T(\o X,\o Y)=\Upsilon\,\varepsilon(\o X)
\varepsilon(\o Y) \o \varepsilon$, where $g(\o \varepsilon,\o
X):=\varepsilon(\o X)$. From which
\begin{eqnarray*}
  S(\o X,\o Y,\o Z,\o W) &=& g(T(\o X,\o W) , T(\o Y, \o Z))-g(T(\o Y, \o W) ,T( \o X, \o
Z)) \\
   &=& \Upsilon\,\varepsilon(\o X)
\varepsilon(\o Y)\varepsilon(\o Z) \varepsilon(\o W)g(\o
\varepsilon,\o \varepsilon)-\Upsilon\,\varepsilon(\o X)
\varepsilon(\o Y)\varepsilon(\o Z) \varepsilon(\o W)g(\o
\varepsilon,\o \varepsilon)=0.
\end{eqnarray*}
\par On the other hand, if the $v$-curvature tensor $S\neq0$, then it
follows from (\ref {8a11}) that $\varepsilon=0$ and $\alpha(\o
\eta)=0$. Hence, $\alpha=0$ and the $hv$-curvature tensor $P$
vanishes. In this case, it follows from the identity (\ref{3.eq.11})
that $\nabla_{\beta \o X}S=0$. \ \ $\Box$

\begin{prop}
A $P_{2}$-like Finsler manifold $(M,L)$ is a $P^*$-Finsler manifold.
\vspace{-0.2cm}
\end{prop}
\prof As $(M,L)$ is
  $P_{2}$-like, then from (\ref{7a5}), we have
$\widehat{P}(\overline{X},\overline{Y})=
  \alpha_{o}T (\overline{X},\overline{Y})$.
Using Lemma \ref{le.p}, we get $(\nabla_{\beta\o\eta}T)(\o X,\o
Y)=\alpha_{0}T(\o X,\o Y)$, from which, by taking the trace,
$\nabla_{\beta\o\eta}C=\alpha_{0}T$, where
$\alpha_{0}=\frac{\widehat{g}(\nabla_{\beta\overline{\eta}}\,C,C)}{C^2}$.
Hence the result. \ \ $\Box$ \vspace{0.2cm}
\par The next definition will be useful in the sequel.
\begin{defn} \vspace{-0.2cm}
A $\pi$-tensor field $\Theta$ is positively homogenous of degree $r$
in  the directional argument $y$ {\em{(}}symbolically,
h{\em{(}}r{\em{))}} if it satisfies the condition
\vspace{-0.2cm}$$\nabla_{\gamma \overline{\eta}}\,\Theta=r\,\Theta,
\ \  or  \ \ D_{\gamma \overline{\eta}}\,\Theta=r\,\Theta.$$
\end{defn}

\begin{lem}\label{hom.}Let $(M,L)$ be a Finsler manifold, then we have
\begin{description}

  \item[(a)] The Finsler metric $g$ {(the angular metric tensor $\hbar$)}
 is homogenous of degree $0$,
        \item[(b)]The $v$-curvature tensor $S$ is homogenous of degree
$-2$,
 \item[(c)]The $hv$-curvature tensor $P$ is homogenous of degree
$-1$,
\item[(d)]The $h$-curvature tensor $R$ is homogenous of degree
$0$,
\item[(e)] The $(h)hv$-torsion tensor $T$ is homogenous of degree
$-1$,
 \item[(f)] The $(v)hv$-torsion tensor $\widehat{P}$ is homogenous of degree
$0$,
 \item[(g)] The $(v)h$-torsion tensor $\widehat{R}$ is homogenous of degree
$1$.
\end{description}
\end{lem}

\begin{lem}\label{nab2} For every vector $(1)\pi$-form $A$, we have\vspace{-0.2cm}
\begin{equation*}
\begin{array}{rcl}
(\stackrel{1}\nabla\stackrel{1}\nabla A)(\o X,\o Y,\o Z)-
(\stackrel{1}\nabla\stackrel{1}\nabla A)(\o Y,\o X,\o Z)&=& A(R(\o
X,\o Y)\o Z)-R(\o X,\o Y)A(\o Z)+\\
&& +(\nabla_{\gamma \widehat{R}(\o X,\o Y)}A)(\o Z).
\end{array}
\end{equation*}
\end{lem}

\vspace{7pt}Deicke theorem \cite{nr2} can be formulated globally as
follows:\vspace{-0.2cm}
\begin{lem} \label{77p} Let  $(M,L)$ be a Finsler manifold. The following
assertions are \linebreak equivalent\,\!{\em:}
\begin{description}
    \item[(a)] $(M,L)$ is Riemannian,
    \item[(b)] The $(h)hv$-torsion tensor $T$ vanishes,
    \item[(c)] The $\pi$-form $C$ vanishes.
\end{description}
\end{lem}

\begin{thm} \label{87p} Let $(M,L)$ be Finsler manifold  which is $h$-isotropic {\em{(}of scalar $k_{0}$)} and
$C^h$-recurrent {\em{(}of recurrence vector $\lambda_{0}$)}. Then,
$(M,L)$ is necessarily one of the following\,\!{\em:}
\begin{description}
    \item[(a)] A Riemannian manifold of constant curvature,
    \item[(b)] A Finsler manifold of dimension 2,
    \item[(c)] A Finsler manifold of dimensions $n\geq3$
     with vanishing scalar $k_{0}$ and\\
     $(\nabla_{\beta \o X}\lambda_{o})(\o Y)=(\nabla_{\beta \o Y}\lambda_{o})(\o
     X)$.
\end{description}
\end{thm}

\prof For a $C^h$-recurrent manifold, one can easily show
that\vspace{-0.1cm}
\begin{eqnarray*}
 & & (\stackrel{1}\nabla\stackrel{1}\nabla T)(\o X,\o Y,\o Z,\o W)
  -(\stackrel{1}\nabla\stackrel{1}\nabla T)(\o Y,\o X,\o Z,\o
  W)=\\&=&
   \{(\nabla_{\beta \o X}\lambda_{o})(\o Y)-(\nabla_{\beta \o Y}\lambda_{o})(\o X)\}T(\o Z,\o
   W)
   =:\Psi(\o X,\o Y)T(\o Z,\o W).\vspace{-0.1cm}
\end{eqnarray*}
From which, taking into account Lemma \ref{nab2}, we obtain
\begin{equation*}\label{74}
  \left.
    \begin{array}{rcl}
  \Psi(\o X,\o Y)T(\o Z,\o W) &=& T(R(\o
X,\o Y)\o Z,\o W)+T(\o Z, R(\o X,\o Y)\o W)-\\
&&-R(\o X,\o Y)T(\o Z,\o W)+(\nabla_{\gamma \widehat{R}(\o X,\o
Y)}T)(\o Z,\o W).\vspace{-0.2cm}
\end{array}
  \right.
\end{equation*}
Now, as $(M,L)$ is $h$-isotropic of scalar $k_{0}$, then the
$h$-curvature tensor $R$ has the form
$$R(\o X,\o Y)\o Z=k_{0}\{g(\o X,\o Z)\o Y-g(\o Y,\o Z)\o X\}\,; \ \ (n\geq3) .  $$
From the above two equations, we get
\begin{equation}\label{71a}
  \left.
   \begin{array}{rcl}
 \Psi(\o X,\o Y)T(\o Z,\o W) &=&k_{0}g(\o X,\o Z)T(\o Y,\o W)-k_{0}g(\o Y,\o Z)T(\o
X,\o W)
 + k_{0}g(\o X,\o W)T(\o Z,\o Y)-\\
 &&- k_{0}g(\o Y,\o W)T(\o Z,\o X)
 -k_{0}g(\o X,T(\o Z,\o W))\o Y+k_{0}g(\o Y,T(\o Z,\o W))\o X\\
&&+ k_{0}g(\o X,\o \eta)(\nabla_{\gamma\o Y}T)(\o Z,\o W)-k_{0}g(\o
Y,\o \eta)(\nabla_{\gamma\o X}T)(\o Z,\o W).\vspace{-0.2cm}
\end{array}
  \right.
\end{equation}
Setting $\o Y=\o \eta$, noting that $T$ is $h(-1)$ and $g(\o \eta
,\o \eta)=L^{2}$, we get
\begin{eqnarray*}
 \Psi(\o X,\o \eta)T(\o Z,\o W) &=&-k_{0}g(\o \eta,\o Z)T(\o
X,\o W)- k_{0}g(\o \eta,\o W)T(\o Z,\o X)
 -k_{0}T(\o X,\o Z,\o W)\o \eta-\\
&&- k_{0}g(\o X,\o \eta)T(\o Z,\o W)-k_{0}L^{2}(\nabla_{\gamma\o
X}T)(\o Z,\o W).\vspace{-0.2cm}
\end{eqnarray*}
From which, we have\vspace{-0.2cm}
\begin{equation}\label{72a}
  \left.
    \begin{array}{rcl}
 g(\o Y, \o \eta)\Psi(\o X,\o \eta)T(\o Z,\o W) &=&-k_{0}g(\o Y, \o \eta)g(\o \eta,\o Z)T(\o
X,\o W)- k_{0}g(\o Y, \o \eta)g(\o \eta,\o W)T(\o Z,\o X)-\\
&& -k_{0}g(\o Y, \o \eta)T(\o X,\o Z,\o W)\o \eta- k_{0}g(\o Y, \o
\eta)g(\o X,\o \eta)T(\o Z,\o
W)-\\
&&-k_{0}L^{2}g(\o Y, \o \eta)(\nabla_{\gamma\o X}T)(\o Z,\o
W),\vspace{-0.2cm}
\end{array}
  \right.
\end{equation}
whereas \vspace{-0.2cm}
\begin{equation}\label{73a}
  \left.
    \begin{array}{rcl}
 g(\o X, \o \eta)\Psi(\o Y,\o \eta)T(\o Z,\o W) &=&-k_{0}g(\o X, \o \eta)g(\o \eta,\o Z)T(\o
Y,\o W)- k_{0}g(\o X, \o \eta)g(\o \eta,\o W)T(\o Z,\o Y)-\\
&& -k_{0}g(\o X, \o \eta)T(\o Y,\o Z,\o W)\o \eta- k_{0}g(\o X, \o
\eta)g(\o Y,\o \eta)T(\o Z,\o
W)-\\
&&-k_{0}L^{2}g(\o X, \o \eta)(\nabla_{\gamma\o Y}T)(\o Z,\o
W).\vspace{-0.2cm}
 \end{array}
  \right.
\end{equation}
Now, from (\ref{71a}), (\ref{72a}) and (\ref{73a}), we obtain
\begin{eqnarray*}
 &&T(\o Z,\o W)\{L^{2}\Psi(\o X,\o Y)-g(\o Y,\o \eta)\Psi(\o X,\o \eta)
 +g(\o X,\o \eta)\Psi(\o Y,\o \eta)\}=\\
 &&
  =\mathfrak{U}_{\o X,\o Y}k_{0}L^{2}\{
  \hbar(\o X,\o Z)T(\o Y,\o W)+\hbar(\o X,\o W)T(\o Y,\o Z)-\phi(\o Y)\,T(\o X,\o Z,\o W)
  \}.
\end{eqnarray*}
Taking the trace of both sides of the above equation, we get
\begin{equation}\label{74a}
  \left.
    \begin{array}{rcl}
&&C(\o Z)\{L^{2}\Psi(\o X,\o Y)-g(\o Y,\o \eta)\Psi(\o X,\o \eta)
 +g(\o X,\o \eta)\Psi(\o Y,\o \eta)\}=\\
 &&
  =2k_{0}L^{2}\{
  \hbar(\o X,\o Z)C(\o Y)- \hbar(\o Y,\o Z)C(\o X)\}.
 \end{array}
  \right.
\end{equation}
Setting $\o Z=\o C$, taking into account the fact that   $\hbar(\o
X,\o C)=C(\o X)$, the above equation reduces to
\begin{equation*}
C(\o C)\{L^{2}\Psi(\o X,\o Y)-g(\o Y,\o \eta)\Psi(\o X,\o \eta)
 +g(\o X,\o \eta)\Psi(\o Y,\o \eta)\}=0.
\end{equation*}
\par Now, if $C(\o C)=g(\o C,\o C)=0$, then $\o C=0$ and so $C=0$.
Consequently, by Lemma \ref{77p}, $(M,L)$ is a Riemannian manifold
of constant curvature.
\par
On the other hand, if  $(M,L)$ is not Riemannian, then we have
\begin{equation*}\label{76a}
L^{2}\Psi(\o X,\o Y)-g(\o Y,\o \eta)\Psi(\o X,\o \eta)
 +g(\o X,\o \eta)\Psi(\o Y,\o \eta)=0.
\end{equation*}
From which, together with (\ref{74a}), we get
\begin{equation}\label{75a}
k_{0}\{\hbar(\o X,\o Z)C(\o Y)- \hbar(\o Y,\o Z)C(\o X)\}=0.
\end{equation}
\par If $k_{0}\neq0$, then, by (\ref{75a}), $\hbar(\o X,\o Z)C(\o Y)=
\hbar(\o Y,\o Z)C(\o X)$. Setting $\o Y=\o C$, we get $\hbar(\o X,\o
Z)=\frac{1}{C^2}C(\o X)C(\o Z)$, which implies that $\dim M=2$.
\par If $k_{0}=0$, then $R=0$ and (\ref{71a}) yields
$\Psi(\o X,\o Y)=0$, which means that $(\nabla_{\beta \o
X}\lambda_{o})(\o Y)=(\nabla_{\beta \o Y}\lambda_{o})(\o X)$. \ \
$\Box$

\vspace{9pt}
Now, we focus our attention  to the interesting case
$(c)$ of the above theorem. In this case, the $h$-curvature tensor
$R=0$ and hence the $(v)h$-torsion tensor $\widehat{R}=0$.
Therefore, the equation (deduced from (\ref{K}))
\begin{equation*}\label{a}
  \left.
    \begin{array}{rcl}
   & &(\nabla_{\gamma\o X}R)(\o Y,\o Z,\o W) + (\nabla_{\beta\o Y}P)(\o
Z,\o X,\o W)-(\nabla_{\beta \o Z}P)(\o Y,\o X,\o W)-\\
& &- P(\o Z,P(\o Y,\o X)\o \eta)\o W+R(T(\o X,\o Y),\o Z)\o W-S(R(\o
Y,\o Z)\o \eta,\o X)\o W+\\ & &+ P(\o Y, P(\o Z,\o X)\o \eta)\o W
-R(T(\o X,\o Z),\o Y)\o W=0. \vspace{-0.1cm}
\end{array}
  \right.
\end{equation*}
reduces to
\begin{equation*}
  \left.
    \begin{array}{rcl}
   & & (\nabla_{\beta\o Y}P)(\o
Z,\o X,\o W)-(\nabla_{\beta \o Z}P)(\o Y,\o X,\o W)-\\ && -P(\o
Z,\widehat{P}(\o Y,\o X))\o W+ P(\o Y, \widehat{P}(\o Z,\o X))\o W
=0. \vspace{-0.1cm}
\end{array}
  \right.
\end{equation*}
Setting $\o W=\o \eta$, we get
\begin{equation}\label{7a1}
  \left.
    \begin{array}{rcl}
   & & (\nabla_{\beta\o Y}\widehat{P})(\o
Z,\o X)-(\nabla_{\beta \o Z}\widehat{P})(\o Y,\o X)- \widehat{P}(\o
Z,\widehat{P}(\o Y,\o X))+ \widehat{P}(\o Y, \widehat{P}(\o Z,\o X))
=0. \vspace{-0.1cm}
\end{array}
  \right.
\end{equation}
Since $(M,L)$  is $C^h$-recurrent, then,  by  Proposition \ref{7p},
 the $(v)hv$-torsion tensor $\widehat{P}$ satisfies the relations  $(
\nabla_{\beta \o Z} \widehat{P})(\o X,\o Y)=(K_{o}\lambda_{o}(\o Z)+
    \nabla_{\beta \o Z} K_{o})T(\o X,\o Y)$ and
    $\widehat{P}(\o X,\o Y)=\lambda_{o}(\o \eta) T(\o X,\o Y)=K_{o}T(\o X,\o
    Y)$. From these, together with  (\ref{7a1}), we get

\begin{equation*}
  \left.
    \begin{array}{rcl}
      & &  (K_{o}\lambda_{o}(\o Y)+
    \nabla_{\beta \o Y} K_{o})T(\o Z,\o X)-(K_{o}\lambda_{o}(\o Z)+
    \nabla_{\beta \o Z} K_{o})T(\o X,\o Y)- \\
    &&-K_{o}^{2}T(\o Z,T(\o X,\o  Y))+ K_{o}^{2}T(\o Y, T(\o X,\o
    Z))=0.
\vspace{-0.1cm}
\end{array}
  \right.
\end{equation*}
Hence, by Corollary \ref{CC},\vspace{-0.1cm}
$$K_{o}^{2}S(\o Y,\o Z,\o X,\o W)=\mathfrak{U}_{\o Y,\o Z}\{(K_{o}\lambda_{o}(\o Y)+
    \nabla_{\beta \o Y} K_{o})T(\o X,\o Z,\o W)\}.\vspace{-0.1cm}$$
\par  As $S(\o Y,\o Z,\o X,\o W)$ is skew-symmetric in the arguments
$\o X$ and $\o W$ while the right-hand side is symmetric in the same
arguments, we obtain
\begin{equation}\label{7a2}
K_{o}^{2}S(\o Y,\o Z,\o X,\o W)=0,\vspace{-0.1cm}
\end{equation}
\begin{equation}\label{7a3}
\mathfrak{U}_{\o Y,\o Z}\{(K_{o}\lambda_{o}(\o Y)+
    \nabla_{\beta \o Y} K_{o})T(\o Z,\o X,\o W)\}=0.\vspace{-0.1cm}
\end{equation}
It follows from (\ref{7a2}) and (\ref{7a4}) that\vspace{-0.1cm}
 $$ P(\o X,\o Y,\o Z, \o W)=\lambda_{o}(\o Z)T(\o X,\o Y,\o W)-\lambda_{o}(\o W)T(\o X,\o Y,\o
    Z)\vspace{-0.1cm}.$$
\par  On the other hand, if $K_{o}\neq0$, then the $v$-curvature tensor
$S$ vanishes from  (\ref{7a2}). Next, it is seen from (\ref{7a3})
that, if $\mathbf{V}(\o Y):=K_{o}\lambda_{o}(\o Y)+
    \nabla_{\beta \o Y} K_{o}\neq0$, then there exists a scalar
    function
     $\Upsilon=\frac {T(\o X,\o Z,\o W)T(\o X,\o Y,\o Z)T(\o Y,\o Z,\o W)}
     {(T(\o X,\o Y,\o W))^{2} (\mathbf{V}(\o Z))^3}$ \  such that
     $$ T(\o X,\o Y, \o W)=\Upsilon\, \mathbf{V}(\o X) \mathbf{V}(\o Y) \mathbf{V}(\o W).$$

 Summing up, we have\vspace{-0.2cm}
 \begin{thm}\label{7th} Let $(M,L)$ be a Finsler manifold of dimensions
 $n\geq3$. If $(M,L)$ is $h$-isotropic  and $C^h$-recurrent, then\vspace{-0.2cm}

 \begin{description}
 \item [(a)] the recurrence vector $\lambda_{o}$ satisfies\;\!\em{:} $(\nabla_{\beta \o X}\lambda_{o})(\o Y)=(\nabla_{\beta \o Y}\lambda_{o})(\o
     X)$,
    \item [(b)] the $h$-curvature tensor $R=0$ and
 the $(v)h$-torsion tensor $\widehat{R}=0$,
    \item [(c)] the $hv$-curvature tensor $P$ has the property that \\
    $ P(\o X,\o Y,\o Z, \o W)=\lambda_{o}(\o
Z)T(\o X,\o Y,\o W)-\lambda_{o}(\o W)T(\o X,\o Y,\o Z),$
    \item [(d)] the $(v)hv$-torsion tensor $\widehat{P}(\o X,\o Y)=K_{o}T(\o X,\o Y)$.
 \end{description}
\vspace{-7pt} \par  Moreover, if $K_{o}\neq0$, then \vspace{-0.2cm}
  \begin{description}
    \item [(e)] the $v$-curvature tensor $S$ vanishes,
    \item [(f)] the $(h)hv$-torsion tensor $T$ satisfies\;\!\em{:} $ T(\o X,\o Y, \o
W)=\Upsilon\, \mathbf{V}(\o X) \mathbf{V}(\o Y) \mathbf{V}(\o W).$
  \end{description}
\end{thm}
By Definition \ref{def.p2like} and Theorem \ref{7th}, we immediately
have\,: \vspace{-0.2cm}
\begin{cor}\label{7c.5} A Finsler manifold $(M,L)$ of dimension
 $n\geq3$ which is $h$-isotropic  and $C^h$-recurrent is necessarily $P_{2}$-like.
\end{cor}

\vspace{0.1cm}
 Now, we define an  operator $\mathbb{P}$ which aids us to investigate the $R_{3}$-like manifolds.\vspace{-0.3cm}
\begin{defn} \label{def.i1} ~\par\vspace{-0.2cm}
\begin{description}
    \item[(a)] If $\omega$ is a $\pi$-tensor field  of type {\em(1,p)},
    then $ \mathbb{P}\cdot \omega$ is a $\pi$-tensor field  of the same type
    defined by{\,\em:}\vspace{-0.3cm}
\begin{equation*}\label{eq.i2}
    (\mathbb{P}\cdot \omega)(\o X_{1},..., \o X_{p}):= \phi(\omega(\phi(\o X_{1}),..., \phi(\o
    X_{p}))),
\end{equation*}
where $\phi$ is the  vector $\pi$-form defined by
{\em(\ref{eq.i1})}.
    \item[(b)] If $\omega$ is a $\pi$-tensor field  of type {\em(0,p)},
    then $ \mathbb{P}\cdot \omega$ is a $\pi$-tensor field  of the same type
    defined by{\,\em:}\vspace{-0.3cm}
\begin{equation*}\label{eq.i2}
    (\mathbb{P}\cdot \omega)(\o X_{1},..., \o X_{p}):= \omega(\phi(\o X_{1}),..., \phi(\o
    X_{p})).
    \end{equation*}
\end{description}
\end{defn}

\begin{rem}\label{3.rem.1} Since $\phi(\phi(\o X))=\phi(\o X)$ for every $\o X\in\cp$ {\em{(Lemma \ref{cor.i1})}},
then the operator $\mathbb{P}$ is a projector  {\em(i.e.
  $\mathbb{P}\cdot(\mathbb{P}\cdot\omega)=\mathbb{P}\cdot\omega$)}.
 \end{rem}

\begin{defn}\label{4.def.1}A $\pi$-tensor field $\omega$ is said to be indicatory
 if it  satisfies the \linebreak condition{\em\,:} $ \mathbb{P}\cdot \omega
     =\omega$.
\end{defn}

The following result gives a characterization of the indicatory
property for certain types of $\pi$-tensor fields\,:\vspace{-0.2cm}
\begin{lem}\label{cor.i2}~\par\vspace{-0.2cm}
\begin{description}
    \item[(a)] A vector {\em(2)}\,$\pi$-form
 $\omega$ is indicatory if, and only if,  $\omega(\o X,\o \eta)=0=\omega(\o \eta, \o X)$
 and $g(\omega(\o X, \o Y),\o \eta)=0 $.
    \item[(b)] A scaler {\em(2)}\,$\pi$-form
  $\omega$ is  indicatory if, and only if,  $\omega(\o X,\o \eta)=0=\omega(\o \eta, \o
  X)$.
\end{description}
\end{lem}

\prof\\
\textbf{(a)} Let $\omega$ be a vector {(2)}\,$\pi$-form. By
Definition \ref{def.i1}(a) and taking into account (\ref{eq.i1}), we
get\vspace{-0.1cm}
\begin{equation}\label{4.eq.1}
  \left.
    \begin{array}{rcl}
       (\mathbb{P}\cdot \omega)(\o X, \o Y) & = &\phi(\omega(\phi(\o X), \phi(\o
       Y)))\\
       &=&\phi\{\omega(\o X-L^{-1}\ell(\o X)\o \eta,\o Y-L^{-1}\ell(\o Y)\o
       \eta)\}\\
       &=&\phi\{\omega(\o X, \o Y) -L^{-1}\ell(\o Y)\omega(\o X, \o
       \eta)-\\
       &&-L^{-1}\ell(\o X)\omega(\o \eta,\o Y)+L^{-2}\ell(\o X)\ell(\o
       Y)\omega(\o \eta,\o \eta)\}\\
       &=&\omega(\o X,\o Y)-L^{-2}g(\omega(\o X,\o Y),\o \eta)\o \eta-
       \phi\{L^{-1}\ell(\o Y)\omega(\o X, \o
       \eta)+\\
       &&+L^{-1}\ell(\o X)\omega(\o \eta,\o Y)-L^{-2}\ell(\o X)\ell(\o
       Y)\omega(\o \eta,\o \eta)\}
       \end{array}
  \right. \vspace{-0.2cm}
\end{equation}
\par Now, if  $\omega(\o X,\o \eta)=0=\omega(\o \eta, \o X)$
 and $g(\omega(\o X, \o Y),\o \eta)=0 $, then (\ref{4.eq.1}) implies that
 $ (\mathbb{P}\cdot \omega)(\o X,\o Y)
     =\omega(\o X,\o Y)$  and hence $\omega$ is
     indicatory.
\par On the other hand, if $\omega$ is indicatory, then
  $\omega(\o X,\o Y)=\phi(\omega(\phi(\o X), \phi(\o
       Y)))$. From which, setting $\o X=\o \eta$ (resp. $\o Y=\o
       \eta$) and taking into account the fact that $\phi(\o \eta)=0$ (Lemma
       \ref{cor.i1}), we get  $\omega(\o \eta,\o Y)=0$ (resp. $\omega(\o X, \o
       \eta)=0)$. From this, together with $ (\mathbb{P}\cdot \omega)(\o X,\o Y)
     =\omega(\o X,\o Y)$, Equation (\ref{4.eq.1}) implies that $L^{-2}g(\omega(\o X,\o Y),\o \eta)\o
     \eta=0$. Consequently, $g(\omega(\o X,\o Y),\o \eta)=0$.

\noindent \textbf{(b)} The proof is similar to that of (a) and we
omit it.\ \ $\Box$

\begin{prop}\label{p.1} For a Finsler manifold $(M,L)$, the  following  tensors
are \linebreak indicatory\,{\em:}\vspace{-0.2cm}
\begin{description}
    \item[(a)] The $\pi$-tensor field $\phi$,
    \item[(b)] The mixed torsion tensor $T$,
    \item[(c)] The $v$-curvature tensor $S$,
    \item[(d)] The angular metric tensor   $\hbar$,
    \item[(e)] The $\pi$-tensor field $\ \mathbb{P}\cdot
    \omega$ for every $\pi$-tensor field $\omega$.
\end{description}
\end{prop}

Now, we define  the following  $\pi$-tensor fields\,\!:
\vspace{-0.2cm}
 \begin{equation}\label{eq.i6}
  \left.
    \begin{array}{rcl}
 F&:& \ \ F(\overline{X},\overline{Y}):={ \frac{1}{n-2}\{
Ric^h(\overline{X},\overline{Y}) - \frac{Sc^h
g(\overline{X},\overline{Y})}{2(n-1)}\}}, \\
 F_{o}&:&\ \ g(F_{o}(\o X), \o Y):= F(\o X,\o Y), \\
F^{a}&:& \ \ F^{a}(\o X):=F(\o\eta, \o X), \\
 F^{b}&:& \ \ F^{b}(\o X):=F(\o X, \o \eta), \\
  m     & :& \ \ m(\o X,\o Y):= (\mathbb{P}\cdot F)(\o X,\o Y),\\
  m_{o}&:&\ \  g(m_{o}(\o X),\o Y):= m(\o X,\o Y),\\
  a &:&  \ \ a(\o X):=L^{-1}(\mathbb{P}\cdot F^{a})(\o X),\\
  \o a&:&\ \ g(\o a,\o Y):= a(\o X),\\
  b  &:& \ \ b(\o X):=L^{-1}(\mathbb{P}\cdot F^{b})(\o X),\\
  \o b&:&\ \ g(\o b,\o X):= b(\o X), \\
  c &:& \ \ c:=L^{-2} F(\o \eta,\o \eta),\\
   \widehat{R}&:& \ \ \widehat{R}(\o X,\o Y):= R(\o X,\o Y)\o \eta ,\\
  H&:&\ \ H(\o X):= R(\o \eta,\o X)\o \eta=\widehat{R}(\o\eta,\o X).
    \end{array}
  \right\}
\end{equation}
\begin{rem}\label{rem.i2} One can show that \;$m$, $m_{o}$, $a$ and \;$b$ are indicatory and $H(\o \eta)=0$.
\end{rem}

\begin{prop}\label{p.2} If  $(M, L)$ is  an  $R_{3}$-like Finsler manifold, then the
$\pi$-tensor field $F$ can be written in the form\vspace{-0.2cm}
\begin{equation}\label{eq.i4}
    F(\o X,\o Y)= m(\o X,\o Y)+\ell(\o X)a(\o Y)+\ell(\o Y)b(\o
    X)+c\,\ell(\o X)\ell(\o Y) .
\end{equation}
\end{prop}

\prof The proof  follows from Definitions \ref{def.r3} and
\ref{def.i1}(b), taking into account Equations (\ref{eq.i1}) and
(\ref{eq.i6}). In more details\,:
\begin{equation*}\label{4.eq.2}
  \left.
    \begin{array}{rcl}
 (\mathbb{P}\cdot F)(\o X,\o Y) &=& F(\phi(\o X), \phi(\o
       Y))\\
       &=&F(\o X-L^{-1}\ell(\o X)\o \eta,\o Y-L^{-1}\ell(\o Y)\o
       \eta)\\
       &=&F(\o X, \o Y) -L^{-1}\ell(\o Y)F(\o X, \o
       \eta)-\\
       &&-L^{-1}\ell(\o X)F(\o \eta,\o Y)+L^{-2}\ell(\o X)\ell(\o
       Y)F(\o \eta,\o \eta)\\
       &=&F(\o X, \o Y) -L^{-1}\ell(\o Y)\{(\mathbb{P}\cdot F^{b})(\o X)
       +L^{-1}\ell(\o X)F(\o \eta ,\o \eta )\}-\\
       &&-L^{-1}\ell(\o X)\{(\mathbb{P}\cdot F^{a})(\o Y)
       +L^{-1}\ell(\o Y)F(\o \eta\,\o \eta )\}
       +L^{-2}\ell(\o X)\ell(\o
       Y)F(\o \eta,\o \eta)\\
       &=&F(\o X,\o Y)-\ell(\o X)a(\o Y)-\ell(\o Y)b(\o
    X)-c\,\ell(\o X)\ell(\o Y).\ \ \Box
    \end{array}
  \right.
\end{equation*}


\begin{rem}\label{rem.i1} One can show that the $\pi$-tensor fields $a$ and $b$
satisfy the following relations\vspace{-0.2cm}
\begin{equation}\label{eq.i5}
  \left.
    \begin{array}{rcl}
F^{a}(\o X) &=& L \{a(\o X)+c\, \ell(\o X)\}, \\
 F^{b}(\o X) &=& L \{b(\o X)+c\, \ell(\o X)\}.
    \end{array}
  \right.
\end{equation}
\end{rem}

\begin{prop}\label{pp.r}In an $R_{3}$-like Finsler manifold $(M,L)$, we
have\,{\em{:}}
\begin{description}
    \item[(a)]$ R(\o X,\o Y)\o Z=g(\o X,\o Z)F_{o}(\o Y)+F(\o X,\o Z) \o Y
                -g(\o Y,\o Z)F_{o}(\o X)-F(\o Y,\o Z)\o X$.

    \item[(b)] $ \widehat{R}(\o X,\o Y)=g(\o X,\o \eta)F_{o}(\o Y)+F(\o X,\o \eta) \o Y
                -g(\o Y,\o \eta)F_{o}(\o X)-F(\o Y,\o \eta)\o X$.
    \item[(c)] $ H(\o Y)=L^{2}F_{o}(\o Y)+c\,L^{2} \o Y
                -g(\o Y,\o \eta)F_{o}(\o \eta)-F(\o Y,\o \eta)\o
                \eta$.

    \item[(d)]$F_{o}(\o X)=m_{o}(\o X)+\o a\, \ell(\o X)+L^{-1}b(\o X)\o \eta
    +c\, L^{-1}\ell(\o X) \o \eta$.
\end{description}
Consequently,
\begin{description}
    \item[(e)] $ \widehat{R}(\o X,\o Y)= L\{\ell(\o X)(m_{o}(\o Y)+c\,\phi(\o Y))+b(\o X) \phi(\o
    Y)\}-$\\
       ${\qquad\qquad\:}  - L\{\ell(\o Y)(m_{o}(\o X)+c\,\phi(\o X))+b(\o Y) \phi(\o
       X)\}$.

    \item[(f)] $ H(\o Y)= L^{2}\{m_{o}(\o Y)+c\,\phi(\o
    Y)\}$.
\end{description}
\end{prop}

\prof\\
\textbf{(a)} Since $(M,L)$ is an $R_{3}$-like manifold, then by
Definition \ref{def.r3}, we have\vspace{-0.2cm}
\begin{equation*}\label{h}
   \begin{split}
    R(\overline{X},\overline{Y},\overline{Z},\overline{W})= &
       g(\overline{X},\overline{Z})F(\overline{Y},\overline{W})
   -g(\overline{Y},\overline{Z})F(\overline{X},\overline{W}) + \\
       &+ g(\overline{Y},\overline{W})F (\overline{X},\overline{Z})
       - g(\overline{X},\overline{W})F(\overline{Y},\overline{Z})
       .
   \end{split}
\end{equation*}
From which, using the fact that $ g(F_{o}(\o X), \o Y)= F(\o X,\o
Y)$ and that the Finsler metric $g$ is non-degenerate, the result
follows. \vspace{5pt}

\noindent\textbf{(b)} Follows from (a) by setting $\o Z=\o \eta$.
\vspace{5pt}

\noindent\textbf{(c)} Follows from (b) by setting $\o X=\o \eta$.

\vspace{5pt} \noindent\textbf{(d)} By (\ref{eq.i4}) and  (\ref{eq.i6}), we get\\
 $g(F_{o}(\o X),\o Y)=g(m_{o}(\o X),\o Y)+g(\o a,\o Y)\, \ell(\o
X)+L^{-1}b(\o X)g(\o \eta,\o Y)
    +c\, L^{-1}\ell(\o X) g(\o \eta,\o Y).$ Hence, the result
    follows, from the non-degeneracy of $g$.

\vspace{5pt} \noindent\textbf{(e)} Follows by substituting $F_{o}(\o
X)$ (from (d)) and $F^{b}(\o X)$ (from (\ref{eq.i5}))
  into  (b).

\vspace{5pt} \noindent\textbf{(f)} Follows from (e) by setting $\o
X=\o \eta$, taking into account Remark \ref{rem.i2} and the fact
that $\ell(\o \eta)=L$. \ \ $\Box$

\begin{rem}
In view of {\em(\ref{eq.i6})} and Lemma {\em\ref{cor.i1}},
Definition
{\em\ref{def.2}(a)} can be reformulated as follows: \\
A Finsler manifold $(M,L)$ is of scaler curvature if  the
$\pi$-tensor field $H$ satisfies the relation \,$H(\o X)=L^{2}\kappa
\, \phi(\o X)$, where $\kappa$ is a scalar function on $\tm.$
\end{rem}

\begin{defn}\label{d.1} A Finsler manifold $(M,L)$ is said to be of perpendicular
scalar \linebreak {\em({\it or of $p$-scalar})} curvature if the
$h$-curvature tensor $R$ satisfies the condition\vspace{-0.2cm}
\begin{equation}\label{eq.i7}
    (\mathbb{P}\cdot R)(\o X, \o Y,\o Z,\o W)=R_{o} \{ \hbar(\o X, \o Z)\hbar(\o Y,\o W)
    - \hbar(\o X,\o W) \hbar(\o Y, \o Z)\} ,\vspace{-0.2cm}
\end{equation}
where $R_{o}$ is a function  called the perpendicular scalar
curvature.
\end{defn}

\begin{defn}\label{def.p2} A Finsler manifold $(M,L)$ is said to be of $s$-$ps$ curvature if
$(M,L)$ is both of scalar curvature and of $p$-scalar curvature.
\end{defn}
\vspace{-0.3cm}
\begin{prop}\label{pp.r2} If $m_{o}(\o X)=t\, \phi(\o X)$, then an $R_{3}$-like Finsler
manifold is a Finsler manifold of $s$-$ps$ curvature.
\end{prop}
\vspace{-0.1cm} \prof Under the given assumption and  taking into
account Proposition \ref{pp.r}(f), we have\vspace{-0.1cm}
\begin{equation*}\label{2.eq1}
    H(\o X)=L^{2}\kappa\phi(\o X) , \ with \ \kappa=t+c .\vspace{-0.1cm}
\end{equation*}
 Thus, the considered manifold is of scalar curvature.
\par
Now, we prove that the given manifold is of $p$-scalar curvature.
Applying the projection $\mathbb{P}$ on the $h$-curvature tensor $R$
of an $R_{3}$-like manifold, we get \vspace{-0.2cm}
  \begin{equation}\label{3.eq.i5}
     \begin{array}{rcl}
 (\mathbb{P}\cdot R)(\o X, \o Y,\o Z,\o W) &=&  R(\phi (\o X), \phi (\o Y),\phi (\o Z),\phi (\o W)) \\
   &=&g(\phi (\o X),\phi (\o Z))(\mathbb{P}\cdot F)(\o Y,\o W)+
   g(\phi (\o Y),\phi (\o W))(\mathbb{P}\cdot F)(\o X,\o Z)-\\
   & &-g(\phi (\o Y),\phi (\o Z))(\mathbb{P}\cdot F)(\o X,\o W)-
   g(\phi (\o X),\phi (\o W))(\mathbb{P}\cdot F)(\o Y,\o Z)\\
    &=&g(\phi (\o X),\phi (\o Z))m(\o Y,\o W)+
   g(\phi (\o Y),\phi (\o W))m(\o X,\o Z)-\\
   & &-g(\phi (\o Y),\phi (\o Z))m(\o X,\o W)-
   g(\phi (\o X),\phi (\o W))m(\o Y,\o Z).
  \end{array}
     \end{equation}
Since\vspace{-0.2cm}
\begin{eqnarray*}
 g(\phi (\o X),\phi (\o Y))&=& g(\phi (\o X),\o Y-L^{-1}\ell(\o Y)\o
 \eta)=g(\phi (\o X),\o Y)-L^{-1}\ell(\o Y)g(\phi (\o X), \o  \eta)\\
  &=&\hbar(\o X,\o Y)-L^{-1}\ell(\o Y)\hbar(\o X, \o  \eta)=\hbar(\o X,\o Y),\vspace{-0.2cm}
\end{eqnarray*}
then, by using again the given assumption ($ m_{o}=t \,\phi
\Longrightarrow m=t \hbar$), Equation (\ref{3.eq.i5}) reduces
to\vspace{-0.1cm}
\begin{eqnarray*}
 (\mathbb{P}\cdot R)(\o X, \o Y,\o Z,\o W)&=& \hbar(\o X, \o Z)m(\o Y,\o W)+
   \hbar(\o Y,\o W)m(\o X,\o Z)-\\
   & &-\hbar(\o Y,\o Z)m(\o X,\o W)-
   \hbar(\o X,\o W)m(\o Y,\o Z)\\
   &=&2t\{ \hbar(\o X, \o Z)\hbar(\o Y,\o W)
    - \hbar(\o Y, \o Z)\hbar(\o X,\o W) \}.\vspace{-0.1cm}
\end{eqnarray*}
Therefore,  by taking $R_{o}=2t$, we have  \vspace{-0.1cm}
\begin{equation*}\label{3.eq.i7}
    (\mathbb{P}\cdot R)(\o X, \o Y,\o Z,\o W)=R_{o}\{ \hbar(\o X, \o Z)\hbar(\o Y,\o W)
    - \hbar(\o Y, \o Z)\hbar(\o X,\o W) \} .\vspace{-0.1cm}
\end{equation*}
Consequently,  the given manifold is of $p$-scalar curvature. \ \
$\Box$

\begin{thm}\label{th.r1}If an $R_{3}$-like Finsler manifold $(M,L)$ is of $p$-scalar
curvature, then it is of $s$-$ps$ curvature.
\end{thm}

\prof Since the considered manifold is $R_{3}$-like, then, by the
same procedure as in the proof of Proposition \ref{pp.r2}, we
have\vspace{-0.1cm}
\begin{equation}\label{eq.i8}
  \left.
    \begin{array}{rcl}
  (\mathbb{P}\cdot R)(\o X, \o Y,\o Z,\o W) &=& \hbar(\o X, \o Z)m(\o Y,\o W)+
   \hbar(\o Y,\o W)m(\o X,\o Z)-\\
   & &-\hbar(\o Y,\o Z)m(\o X,\o W)-
   \hbar(\o X,\o W)m(\o Y,\o Z).
    \end{array}
  \right.
\end{equation}
On the other hand, since the considered manifold is of $p$-scalar
curvature, then the $h$-curvature tensor satisfies \vspace{-0.2cm}
\begin{equation}\label{eq.i9}
    (\mathbb{P}\cdot R)(\o X, \o Y,\o Z,\o W)=R_{o}\{ \hbar(\o X, \o Z)\hbar(\o Y,\o W)
    - \hbar(\o Y, \o Z)\hbar(\o X,\o W) \}. \vspace{-0.2cm}
\end{equation}
Now, from Equations (\ref{eq.i8}) and (\ref{eq.i9}), we obtain
\vspace{-0.2cm}
\begin{equation*}
   \mathfrak{U}_{\o X,\o Y}\{ R_{o} \hbar(\o X, \o Z)\hbar(\o Y,\o W)-
    \hbar(\o X, \o Z)m(\o Y,\o W)-\hbar(\o Y,\o W)m(\o X,\o Z)\}=0. \vspace{-0.2cm}
\end{equation*}
Using (\ref{eq.i6}) and the non-degeneracy of the metric tensor $g$,
the above equation reduces to \vspace{-0.2cm}
\begin{equation}\label{eq.i10}
   \mathfrak{U}_{\o X,\o Y}\{ R_{o} \hbar(\o X, \o Z)\phi(\o Y)-
    \hbar(\o X, \o Z)m_{o}(\o Y)-m(\o X,\o Z)\phi(\o Y)\}=0. \vspace{-0.2cm}
\end{equation}
Since the $\pi$-tensor fields $\phi, m$ and $m_{o}$ are indicatory,
then \\
$Tr\{\o Y\longmapsto \hbar(\o X, \o Y)\phi(\o Z)\}=g(\o X,\phi(\o Z)
)=\hbar(\o X,\o Z)$,\\
$Tr\{\o Y\longmapsto \hbar(\o X, \o Y)m_{o}(\o Z)\}=m(\o X,\o Z)$,\\
$Tr\{\o Y\longmapsto m(\o X, \o Y)\phi(\o Z)\}=m(\o X,\o Z)$.\\
Consequently, if we take the trace of both sides of Equation
(\ref{eq.i10}), making use of Lemma \ref{cor.i2}, we
get\vspace{-0.2cm}
\begin{equation*}
(n-2)R_{o}\hbar(\o X, \o Z)-(n-3)m(\o X, \o Z)-(n-1)t\,\hbar(\o X,
\o Z)=0, \vspace{-0.2cm}
\end{equation*}
where \,$t:=\frac{1}{n-1}\, Tr( m_{o})$. From which, using
(\ref{eq.i6}) and Lemma \ref{cor.i1}, we get \vspace{-0.2cm}
\begin{equation}\label{eq.i11}
(n-2)R_{o}\phi-(n-3)m_{o}-(n-1)t\,\phi=0. \vspace{-0.2cm}
\end{equation}
Again, taking the trace of the above equation, we obtain
\vspace{-0.2cm}$$(n-1)(n-2) (R_{o}-2t)=0.\vspace{-0.2cm}$$
Substituting the above relation into (\ref{eq.i11}), we get
$\,m_{o}=t\,\phi$. Hence, by  Proposition \ref{pp.r2}, the result
follows. \ \ $\Box$

\begin{thm}\label{th.7}If an $R_{3}$-like Finsler manifold $(M,L)$ is of scalar
curvature, then it is of $s$-$ps$ curvature.
\end{thm}

\prof Since the given manifold is $R_{3}$-like, then the
$\pi$-tensor $H$ is given by (cf. Proposition \ref{pp.r}):
\begin{equation}\label{eq.i13}
H(\o X)= L^{2}\{m_{o}(\o X)+c\,\phi(\o X)\}.
\end{equation}
And since the considered manifold is of scalar curvature, then
\begin{equation}\label{eq.i14}
H(\o X)= L^{2}\kappa \phi(\o X).
\end{equation}
From Equations (\ref{eq.i13}) and (\ref{eq.i14}), we deduce that
$m_{o}(\o X)=(\kappa-c)\phi(\o X)=:t \phi(\o X)$. \\
Hence, by Proposition \ref{pp.r2}, the result follows.\ \ $\Box$

\vspace{8pt}
 Now, let us define the $\pi$-tensor field
\begin{equation}\label{eq.i12}
  \left.
\begin{array}{rcl}
\Psi(\o X, \o Y,\o Z,\o W)&=&R(\o X, \o Y,\o Z,\o W)-
 \frac{1}{n-2}\mathfrak{U}_{\o X, \o Y}\{ g(\o X, \o Z)Ric^{h}(\o Y,\o W)+\\
   & &+g(\o Y,\o W)Ric^{h}(\o X,\o Z)-r g(\o X,\o Z)g(\o Y,\o W)\},
\end{array}
  \right.
\end{equation}
 \vspace{-0.3cm}
where $r=\frac{1}{n-1}\,Sc^{h}$. From Definition \ref{def.r3} and
(\ref{eq.i12}), we immediately obtain

\begin{thm}\label{th.6}
An $R_{3}$-like Finsler manifold is characterized by
$$\Psi(\o X,\o Y,\o Z, \o W)=0.$$
\end{thm}
The tensor field $\Psi$ in the above theorem being of the same form
as the Weyl conformal tensor in  Riemannian geometry, we draw the
following
\begin{thm}\label{th.6a}
An $R_{3}$-like Riemannian manifold is conformally flat.
\end{thm}
\begin{rem}It should be noted that some important results of
  {\em ~\cite{r75}, \cite{r84}, \cite{r69}, \cite{r65}, \cite{r31},
  \cite{r68}},...,etc.
{\em{(}}obtained in local coordinates{\em{)}} are retrieved from the
above mentioned global results {\em{(}}when localized{\em{)}}.
\end{rem}

\newpage
\vspace{1cm}
\begin{center}
{\bf\Large{\textbf{Appendix. Local formulae}}}\end{center}

 \vspace{3pt}
 \par For the sake of completeness, we present in this appendix
 a brief and concise survey of the local expressions of some
 important geometric objects and the local definitions of the special Finsler manifolds treated in the paper.
\vspace{5pt}
\par
  Let $(U,(x^{i}))$ be  a system  of local coordinates on
 $M$ and $(\pi^{-1}(U),(x^i,y^i))$ the associated system of local coordinates on $TM$.
 We use the following notations\,:\\
 $(\pa_{i}):=(\frac{\pa}{\pa x^i})$: the natural basis of $T_{x}M,\, x\in
 M$,\\
 $(\paa_{i}):=(\frac{\pa}{\pa y^i})$: the natural basis of $V_{u}(\T M),\, u\in
 \T M$,\\
$(\pa_{i},\paa_{i})$: the natural basis of $T_{u}(\T M)$,\\
$(\o \pa_{i} )$: the natural basis of the fiber over $u$ in $\p$
($\o \pa_{i} $ is the lift of $\pa_{i}$ at $u$).
\vspace{5pt}
\par To  a Finsler manifold $(M,L)$, we associate the geometric
objects\,:\\
$g_{ij}:= \frac{1}{2}\,\paa_{i} \paa_{j}L^{2}= \paa_{i} \paa_{j}E$:
the Finsler metric tensor,\\
$C_{ijk}:= \frac{1}{2}\, \paa_{k}\,g_{ij}$: the Cartan tensor,\\
$\hbar_{ij}:= g_{ij}-\ell_{i} \ell_{j}$  \,($\ell_{i}:=\partial
L/\partial y^{i})$: the angular metric tensor,\\
$G^h$: the components of the canonical spray,\\
$G^{h}_{i}:=\paa_{i}G^h$,\\
$G^{h}_{ij}:=\paa_{j}G^h_{i}=\paa_{j}\paa_{i}G^h$,\\
$(\delta_{i}):=(\pa_{i}-G^{h}_{i}\paa_{h})$: the basis of $H_{u}(\T
M)$ adapted to $G^{h}_{i}$,\\
$(\delta_{i}, \paa_{i})$: the basis of $T_{u}(\T M)=H_{u}(\T
M)\oplus V_{u}(\T M)$ adapted to $G^{h}_{i}$.
\vspace{5pt}
\par We have\,:\\
$\gamma(\o \pa_{i})=\paa_{i}$,\\
$\rho(\pa_{i})=\o \pa_{i}$, \ \ $\rho(\paa_{i})=0$,\ \
$\rho(\delta_{i})=\o \pa_{i}$,\\
$\beta(\o \pa_{i})=\delta_{i}$,\\
$J(\pa_{i})= \paa_{i}$, \ \ $J(\paa_{i})=0$,\ \
$J(\delta_{i})= \paa_{i}$,\\
$h:= \beta o \rho= dx^{i} \otimes \pa_{i}- G^{i}_{j}\, dx^{j}
\otimes \paa_{i}$\ \    $v:=\gamma o K=dy^{i} \otimes \paa_{i}+
G^{i}_{j}\, dx^{j} \otimes \paa_{i} $.
\vspace{5pt}
\par We define\,:\\
\noindent $\gamma^{h}_{ij}:= \frac{1}{2}\,g^{h\ell}(\pa_{i}\,g_{\ell
j}+\pa_{j}\,g_{i\ell }- \pa_{\ell}\,g_{i j}
 )$,\\
  $C^{h}_{ij}:= \frac{1}{2}\,g^{h\ell}(\paa_{i}\,g_{\ell j}+\paa_{j}\,g_{i\ell }-
   \paa_{\ell}\,g_{i j})=
  \frac{1}{2}\,g^{h\ell}\,\paa_{i}\,g_{j\ell}=g^{h\ell}\,C_{ij\ell}$,\\
  $\Gamma^{h}_{ij}:= \frac{1}{2}\,g^{h\ell}(\delta_{i}\,g_{\ell j}+\delta_{j}\,g_{i\ell }-
  \delta_{\ell}\,g_{i j})$ .
\vspace{5pt}
  \par
   Then, we have\,:\\
   $\bullet$ The canonical spray $G$: $G^h=\frac{1}{2}\,\gamma^{h}_{ij}\,y^i
   y^j$.\\
  $\bullet$ The Barthel connection
    $\Gamma$: $G^{h}_{i}=
    \paa_{i}G^{h}=\Gamma^{h}_{ij}y^{j}=G^{h}_{ij}y^{j}$.\\
 $\bullet$ The Cartan connection  $C\Gamma$:
  $(\ \Gamma^{h}_{ij}, G^{h}_{i}, \ C^{h}_{ij})$.\\
  The associated $h$-covariant (resp. $v$-covariant) derivative is denoted by
   $\shortmid$ (resp. $\mid$),
   where $K^{i}_{j|k}:=\delta_{k}K^{i}_{j}+K^{m}_{j}\Gamma^{i}_{mk}
   -K^{i}_{m}\Gamma^{m}_{jk}$ \,
   and \,$K^{i}_{j}|_{k}:=\paa_{k}K^{i}_{j}+K^{m}_{j}C^{i}_{mk}
   -K^{i}_{m}C^{m}_{jk}$.\\
$\bullet$ The Berwald connection  $B\Gamma$: $(\ G^{h}_{ij},
G^{h}_{i},\ 0)$.\\
 The associated $h$-covariant (resp.
$v$-covariant) derivative is denoted by
$\stackrel{*}{\shortmid}$(resp. $\stackrel{*}{\mid}$), where
$K^{i}_{j\stackrel{*}{\shortmid}k}:=\delta_{k}K^{i}_{j}+K^{m}_{j}G^{i}_{mk}-K^{i}_{m}G^{m}_{jk}$
   \, and \,$K^{i}_{j}{\stackrel{*}{\mid}}_{k}:=\paa_{k}K^{i}_{j}$.
 \par
 We also have $G^{h}_{ij}= \Gamma^{h}_{ij}+
C^{h}_{ij}{_{|k}}\,y^k= \Gamma^{h}_{ij}+ C^{h}_{ij}{_{|o}}$,\, where
$ C^{h}_{ij}{_{|o}}= C^{h}_{ij}{_{|k}}\,y^{k}$.

\vspace{5pt}
\par
For the Cartan connection, we have\,: \\
$(v)h$-torsion\,:
$R^{i}_{jk}=\delta_{k}G^{i}_{j}-\delta_{j}G^{i}_{k}=\mathfrak{U}_{jk}\{
\delta_{k}G^{i}_{j}\}$,\\
$(v)hv$-torsion\,: $P^{i}_{jk}=
G^{i}_{jk}-\Gamma^{i}_{jk}=C^{i}_{jk|m}y^{m}=C^{i}_{jk|0}$,\\
$(h)hv$-torsion\,: $C^{i}_{jk}={1}/{2}\{ g^{ri}\paa_{r}g_{jk}\}$,\\
$h$-curvature\,:  $R^{i}_{hjk}=\mathfrak{U}_{jk}\{
\delta_{k}\Gamma^{i}_{hj}+\Gamma^{m}_{hj}
\Gamma^{i}_{mk}\}-C^{i}_{hm}
R^{m}_{jk}$,\\
$hv$-curvature\,:
$P^{i}_{hjk}=\paa_{k}\Gamma^{i}_{hj}-C^{i}_{hk|j}+C^{i}_{hm}P^{m}_{jk}$,\\
$v$-curvature\,:
$S^{i}_{hjk}=C^{m}_{hk}C^{i}_{mj}-C^{m}_{hj}C^{i}_{mk}=\mathfrak{U}_{jk}\{
C^{m}_{hk}C^{i}_{mj}\}$.

 \vspace{5pt}
\par
For the Berwald connection, we have\,: \\
$(v)h$-torsion\,:
${R^{*}}^{i}_{jk}=\delta_{k}G^{i}_{j}-\delta_{j}G^{i}_{k}=\mathfrak{U}_{jk}\{
\delta_{k}G^{i}_{j}\}$,\\
$h$-curvature\,:  ${R^{*}}^{i}_{hjk}=\mathfrak{U}_{jk}\{
\delta_{k}G^{i}_{hj}+G^{m}_{hj} G^{i}_{mk}\}$,\\
$hv$-curvature\,:
${P^{*}}^{i}_{hjk}=\paa_{k}G^{i}_{hj}=:G^{i}_{hjk}$.

\vspace{9pt}
 In the following, we give the {\it\textbf{local}} definitions of the special
  Finsler spaces treated in the paper. For each special Finsler space $(M,L)$,
   we set its name, its defining property
   and a selected reference in which the local definition is located:
\begin{itemize}
  \item Rimaniann manifold {\cite{r42}}\,:
  $g_{ij}(x,y)\equiv g_{ij}(x)$ $\Longleftrightarrow$ $C_{ijk}=0$ \
   $\Longleftrightarrow$ $C_{i}:=C^{k}_{ik}=0$ (Deicke's theorem {\cite{nr2}}).

   \item  Minkowaskian manifold {\cite{r42}}:
  $g_{ij}(x,y)\equiv g_{ij}(y)$ $\Longleftrightarrow$ $C^{i}_{jk|h}=0$ \ and\
  ${R}_{ijk}^{h}=0$.

  \item Berwald manifold {\cite{r42}}:
  $\Gamma^{h}_{ij}(x,y)\equiv \Gamma^{h}_{ij}(x)$\,(i.e.\,$\paa_{k}\Gamma^{h}_{ij}=0$)
    $\Longleftrightarrow$ $C^{h}_{ij|k}=0$.

    \item $C^{h}$-recurrent manifold {\cite{r65}}:
     $C_{hij|k}={\mu}_{k}C_{hij}$,\\
      where
    $\mu_{j}$ is a covariant vector field.

    \item $P^{*}$-Finsler manifold {\cite{r76}}:
    $C^{h}_{ij|0}=\lambda(x,y)C^{h}_{ij},$\\
     where
  $\lambda(x,y)=\frac{P_{i}C^{i}}{C^{2}}$;\,
  $P_{i}:=P^{k}_{ik}=C^{k}_{ik|0}=C_{i|0}$ and $C^{2}=C_{i}C^{i}\neq0$.

\item $C^{v}$-recurrent manifold {\cite{r65}}: $C^{i}_{jk}|_{l}=\lambda_{l}\,C^{i}_{jk} \,$ or $\, C_{ijk}|_{l}=
\lambda_{l}\,C_{ijk}$.

\item $C^{0}$-recurrent manifold {\cite{r65}}:
$C^{i}_{jk}\!\!\stackrel{*}{|}_{l}=\lambda_{l}\,C^{i}_{jk}\, $ or
$\, C_{ijk}\!\!\stackrel{*}{\mid}_{l}=\lambda_{l}\,C_{ijk}$.

   \item Semi-$C$-reducible manifold $(\dim M \geq 3)$ {\cite{r32}}:
   \vspace{-0.2cm}$$C_{ijk}= \frac{\mu}{(n+1)}(\hbar_{ij}C_{k}
   +\hbar_{jk}C_{i}+\hbar_{ki}C_{j})+
   \frac{\tau}{C^{2}}C_{i} C_{j} C_{k}, \,\, C^{2}\neq0,$$\vspace{-0.2cm}
    where $\mu$ and $\tau$ are scalar functions satisfying
$\mu +\tau=1$.

    \item $C$-reducible manifold $(\dim M \geq 3)$ {\cite{r2}}:
 $C_{ijk}= \frac{1}{n+1}(\hbar_{ij}C_{k}
  +\hbar_{jk}C_{i}+\hbar_{ki}C_{j}).$

\item $C_{2}$-like manifold $(\dim M \geq 2)$ {\cite{nr3}}:
  $C_{ijk}=  \frac{1}{C^{2}}C_{i} C_{j} C_{k}, \,\, C^{2}\neq0$.

\item quasi-$C$-reducible manifold $(\dim M \geq 3)$ {\cite{nr4}}:
 $C_{ijk}=A_{ij}C_{k}+ A_{jk}C_{i}+ A_{ki}C_{j}, $ \\
 where $A_{ij}(x,y)$
is a symmetric tensor field satisfying  $A_{ij}y^{i}=0$.

\item  $S_{3}$-like manifold $(\dim M\geq 4)$ {\cite{r10}}:
 $S_{lijk}= \frac{S}{(n-1)(n-2)}
\{\hbar_{ik} \hbar_{lj}-\hbar_{ij}\hbar_{lk} \},$\\
where $S$ is the vertical scalar curvature.


\item $S_{4}$-like manifold $(\dim M\geq 5)$ {\cite{r10}}: $
S_{lijk}=\hbar_{lj}\textbf{F}_{ik}-\hbar_{lk}\textbf{F}_{ij}+\hbar_{ik}\textbf{F}_{lj}-\hbar_{ij}\textbf{F}_{lk},$\\
where $\textbf{F}_{ij}:=\frac{1}{n-3} \{S_{ij}-\frac{1}{2(n-2)} S
\hbar_{ij}\}$; $S_{ij}$ being the vertical Ricci tensor.

\item  $S^{v}$-recurrent manifold {\cite{r68},~\cite{r69}}: $S_{hijk}|_{m}=\lambda_{m}
S_{hijk},$\\
 where $\lambda_{j}(x,y)$ is a covariant vector field.

\item  Second order $S^{v}$-recurrent manifold {\cite{r68},~\cite{r69}}:
 $S_{hijk}|_{m}|_{n}=\Theta_{mn}
S_{hijk},$ \\
where $\Theta_{ij}(x,y)$ is a covariant tensor field.

\item Landsberg manifold {\cite{r76}}:
$P^{h}_{kji}\,y^{k}=0\Longleftrightarrow(\dot{\partial_{i}}\Gamma^{h}_{jk})y^{k}=0$
  $\Longleftrightarrow$  $C^{h}_{ij|k}y^{k}=0$.

 \item General Landsberg manifold {\cite{r9}}: $P^{r}_{ijr}y^{i}=0$\
$\Longleftrightarrow$ $C_{j|o}=0$.

 \item $P$-symmetric manifold {\cite{r31}}: $P_{hijk}=P_{hikj}$.

 \item $P_{2}$-like manifold $(\dim M \geq 3)$ {\cite{r29}}:
$P_{hijk}=\alpha_{h}C_{ijk} -\alpha_{i} C_{hjk},$\\
 where  $\alpha_{k}(x,y)$ is a covariant vector  field.
   \item $P$-reducible manifold $(\dim M \geq 3)$ {\cite{r31}}:
$P_{ijk}=\frac{1}{n+1}(\hbar_{ij}\,P_{k}+ \hbar_{jk}\,P_{i}+ \hbar_{ki}\,P_{j}),$\\
 where  $P_{ijk}=g_{hi}P^{h}_{jk}$.

\item $h$-isotropic manifold $(\dim M \geq 3)$ {\cite{r65}}:
 $R_{hijk}= k_{o} \{g_{hj} g_{ik} - g_{hk} g_{ij} \}$,\\
  for some scalar $k_{o}$, where $R_{hijk}=g_{il}R^{l}_{hjk}$.

 \item Manifold of scalar curvature {\cite{nr1}}:
   $R_{ijkl}\,y^{i}y^{k}= k L^{2} \hbar_{jl}$,\\ for some function $k:
\T M \To
 \Real$ .

\item  Manifold of constant curvature {\cite{nr1}}:  the function $k$
in the above definition  is constant.

\item  Manifold of perpendicular scalar {({\it or of $p$-scalar} )} curvature
{\cite{r75},~\cite{r84}}: \\
$\mathbb{P}\cdot R_{hijk}:=\hbar^{l}_{h} \, \hbar^{m}_{i} \,
\hbar^{n}_{j} \, \hbar^{r}_{k} \,R_{lmnr} = R_{o}\{\hbar_{ik}
\hbar_{hj}-\hbar_{ij}\hbar_{hk} \},$\\
 where $R_{o}$ is a function
called a perpendicular scalar curvature.

\item  Manifold of $s$-$ps$ curvature {\cite{r75},~\cite{r84}}:
 $(M,L)$ is both of scalar curvature and of $p$-scalar curvature.

 \item $R_{3}$-like manifold $(\dim M\geq4)$ {\cite{r75}}:
$R_{hijk}=g_{hj}F_{ik}-g_{hk}F_{ij}+g_{ik}F_{hj}-g_{ij}F_{hk},$
\\
where $F_{ij}:= \frac{1}{n-2}\{ R_{ij}-\frac{1}{2}\,r\,g_{ij}\}$;
$R_{ij}:=R^{h}_{ijh}$, \,$r:=\frac{1}{n-1}R^{i}_{i}$.

\end{itemize}

\newpage
\providecommand{\bysame}{\leavevmode\hbox
to3em{\hrulefill}\thinspace}
\providecommand{\MR}{\relax\ifhmode\unskip\space\fi MR }
\providecommand{\MRhref}[2]{%
  \href{http://www.ams.org/mathscinet-getitem?mr=#1}{#2}
} \providecommand{\href}[2]{#2}

\end{document}